\documentclass[12pt]{article}

\DeclareMathAlphabet{\mathpzc}{OT1}{pzc}{m}{it}

\usepackage{amsmath}
\usepackage{amssymb}
\usepackage{amsfonts}
\usepackage{amsthm}
\usepackage{mathrsfs}
\usepackage{ifsym}
\usepackage{fontenc}
\usepackage{textcomp}
\usepackage{tipa}
\usepackage{enumerate}
\usepackage[pdftex]{color, graphicx}

\begin{document}

\title{{\bf  Solutions of diophantine equations as periodic points of $p$-adic algebraic functions, II: The Rogers-Ramanujan continued fraction}}         
\author{Patrick Morton}        
\date{July 16, 2018}          
\maketitle

\begin{abstract}  In this part we show that the diophantine equation $X^5+Y^5=\varepsilon^5(1-X^5Y^5)$, where $\varepsilon=\frac{-1+\sqrt{5}}{2}$, has solutions in specific abelian extensions of quadratic fields $K=\mathbb{Q}(\sqrt{-d})$ in which $-d \equiv \pm 1$ (mod $5$).  The coordinates of these solutions are values of the Rogers-Ramanujan continued fraction $r(\tau)$, and are shown to be periodic points of an algebraic function.
\end{abstract}

\section{Introduction.}

In a previous paper \cite{lym} integral solutions of the diophantine equation
$$Fer_4: \ 16X^4+16Y^4=X^4Y^4,$$
were constructed in ring class fields $\Omega_f$ of odd conductor $f$ over imaginary quadratic fields of the form $K=\mathbb{Q}(\sqrt{-d})$, with $d_Kf^2=-d \equiv 1$ (mod $8$), where $d_K$ is the discriminant of $K$.  The coordinates of these solutions were studied in Part I of this paper \cite{mor5}, and shown to be the periodic points of a fixed $2$-adic algebraic function on the maximal unramified algebraic extension $\textsf{K}_2$ of the $2$-adic field $\mathbb{Q}_2$.  In particular, every ring class field of odd conductor over $K=\mathbb{Q}(\sqrt{-d})$ with $-d \equiv 1$ (mod $8$) is generated over $\mathbb{Q}$ by some periodic point of this algebraic function.  This was simplified and extended in \cite{mor7} to show that all ring class fields over any field $K$ in this family of quadratic fields are generated by individual periodic or pre-periodic points of the $2$-adic algebraic function
$$F(z)=\frac{-1+\sqrt{1-z^4}}{z^2}=\sum_{n=1}^\infty{(-1)^n{\frac{1}{2} \atopwithdelims ( ) n} z^{4n-2}}.$$
A similar situation holds for the solutions of
$$Fer_3: \ 27X^3+27Y^3=X^3Y^3,$$
studied in \cite{mor4}, in that they are, up to a finite set, the exact set of periodic points of a fixed $3$-adic algebraic function, and all ring class fields of quadratic fields $K=\mathbb{Q}(\sqrt{-d})$ in the family for which $-d \equiv 1$ (mod $3$) are generated by periodic or pre-periodic points of this same $3$-adic algebraic function.  (See \cite{mor4} and \cite{mor7} for a more precise description.)  \medskip

In this paper I will study the analogous quintic equation
$$\mathcal{C}_5: \ \upsilon^5X^5+\upsilon^5Y^5=1-X^5Y^5, \quad \upsilon=\frac{1+\sqrt{5}}{2},$$
which can be written in the equivalent form
$$\mathcal{C}_5: \ X^5+Y^5=\varepsilon^5(1-X^5Y^5), \quad \varepsilon=\frac{-1+\sqrt{5}}{2},\eqno{(1.1)}$$
in certain abelian extensions of imaginary quadratic fields $K=\mathbb{Q}(\sqrt{-d})$ with $d_Kf^2=-d \equiv \pm 1$ (mod $5$).  In Part I \cite{mor5} these were called {\it admissible quadratic fields} for the prime $p=5$: these are the imaginary quadratic fields in which the ideal $(5)=\wp_5 \wp_5'$ of the ring of integers $R_K$ of $K$ splits into two distinct prime ideals.  In this part I will show that (1.1) has unit solutions in the abelian extensions $\Sigma_{5} \Omega_f$ or $\Sigma_{5} \Omega_{5f}$ of $K$ (according as $d \neq 4f^2$ or $d=4f^2>4$), where $\Sigma_5$ is the {\it ray class field} of conductor $\mathfrak{f}=(5)$ over $K$ and $\Omega_f, \Omega_{5f}$ are the {\it ring class fields} of conductors $f$ and $5f$, respectively, over $K$, for any positive integer $f$ which is relatively prime to $p=5$.  (See  \cite{co}.)  \medskip

As in the quadratic families mentioned above, the coordinates of these solutions will be shown in Part III to be the exact set of periodic points (minus a finite set) of a specific $5$-adic algebraic function in a suitable extension of the  $5$-adic field $\mathbb{Q}_5$.  This will be used to verify the conjectures of Part I of this paper for the prime $p=5$, according to which any ring class field of conductor $f$ over $K$ with $(f,5)=1$ is generated over $K$ by a periodic point of a fixed $5$-adic algebraic function, independent of $d$, and any ring class field whose conductor is {\it divisible} by $5$ is generated over $K$ by a pre-periodic point of the same algebraic function.  \medskip

Let $H_{-d}(x)$ be the class equation for a discriminant $-d \equiv \pm 1$ (mod $5$), and let
$$F_d(x)=x^{5h(-d)}(1-11x-x^2)^{h(-d)} H_{-d}(j_5(x)),\eqno{(1.2)}$$
where
$$j_5(b)=\frac{(1-12b+14b^2+12b^3+b^4)^3}{b^5(1-11b-b^2)}. \eqno{(1.3)}$$ \smallskip
This rational function represents the $j$-invariant of the Tate normal form
$$E_5(b): \ Y^2 + (1+b)XY +bY=X^3+bX^2,\eqno{(1.4)}$$ \smallskip
on which the point $P=(0,0)$ has order $5$.  Note that
$$j_5(b)=-\frac{(z^2+12z+16)^3}{z+11}, \ \ z=b-\frac{1}{b}.\eqno{(1.5)}$$ \smallskip

The roots of $F_d(x)$ are the values of $b$ for which the curve $E_5(b)$ has complex multiplication by the order $\textsf{R}_{-d}$ of discriminant $-d = d_K f^2$ in $K$.  If $h(-d)$ is the class number of $\textsf{R}_{-d}$, it turns out that $F_d(x^5)$ has an irreducible factor $p_d(x)$ of degree $4h(-d)$ whose roots give solutions of $\mathcal{C}_5$ in abelian extensions of $K=\mathbb{Q}(\sqrt{-d})$.  Furthermore, the roots of $p_d(x)$ are conjugate values over $\mathbb{Q}$ of the Rogers-Ramanujan continued fraction $r(\tau)$ defined by
\begin{align*}
r(\tau) &= \frac{q^{1/5}}{1+\frac{q}{1+\frac{q^2}{1+ \frac{q^3}{1+\cdots}}}}=\frac{q^{1/5}}{1+} \ \frac{q}{1+} \ \frac{q^2}{1+} \ \frac{q^3}{1+} \dots,\\
& = q^{1/5} \prod_{n \ge 1}{(1-q^n)^{(n/5)}},  \ \ q = e^{2 \pi i \tau}, \ \ \tau \in \mathbb{H}.
\end{align*}
See  \cite{anb},  \cite{ber},  \cite{bcz},  \cite{du}.  (We follow the notation in  \cite{du}.)  In the latter formula $(n/5)$ is the Legendre symbol and $\mathbb{H}$ denotes the upper half-plane.  The function $r(\tau)$ is a modular function for the congruence group $\Gamma(5)$ (\cite{du}, p. 149), and $(X,Y)=(r(\tau/5), r(-1/\tau))$ is a modular parametrization of the curve $\mathcal{C}_5$ (see \cite{du}, eq. (7.3)).  In Section 4 we prove the following result. \bigskip

\noindent {\bf Theorem 1.1.}  {\it Let $d \equiv \pm 1$ (mod $5$), $K=\mathbb{Q}(\sqrt{-d})$, and
$$w=\frac{v+\sqrt{-d}}{2} \in R_K, \ \textrm{with} \ \wp_5^2 \mid w \ \textrm{and} \ (N(w),f)=1.$$
Then the values $X=r(w/5), Y=r(-1/w)$ of the Rogers-Ramanujan continued fraction give a solution of $\mathcal{C}_5$ in $\Sigma_{5} \Omega_f$ or $\Sigma_{5} \Omega_{5f}$, according as $d \neq 4f^2$ or $d=4f^2$.  For a unique primitive $5$-th root of unity $\zeta^j=e^{2 \pi i j/5}$, depending on $w$, we have
$$\mathbb{Q}(r(w/5)) = \Sigma_{\wp_5'} \Omega_f, \ \ \mathbb{Q}(\zeta^j r(-1/w)) = \Sigma_{\wp_5} \Omega_f, \ \ \textrm{if} \ d \neq 4f^2;$$
and
$$\mathbb{Q}(r(w/5)) = \Sigma_{2\wp_5'} \Omega_f, \ \ \mathbb{Q}(\zeta^j r(-1/w)) = \Sigma_{2\wp_5} \Omega_f, \ \ \textrm{if} \ d=4f^2, \ 2 \mid f;$$
where $\wp_5$ is the prime ideal $\wp_5 = (5, w)$, $\wp_5'$ is its conjugate ideal in $K$, and $\Sigma_\mathfrak{f}$ denotes the ray class field of conductor $\mathfrak{f}$ over $K$.  Furthermore,
$$\mathbb{Q}(r(-1/w)) = \mathbb{Q}(r(w)) = \Sigma_{5} \Omega_f \ \ \textrm{or} \ \ \Sigma_{5} \Omega_{5f},$$
according as $d \neq 4f^2$ or $d=4f^2$.}
\bigskip

The numbers $\eta=r(w/5), \xi=\zeta^j r(-1/w)$ in this theorem are both roots of the irreducible polynomial $p_d(x)$, and so are conjugate algebraic integers (and units) over $\mathbb{Q}$.  Furthermore, they satisfy the relation
$$\xi =\zeta^j r(-1/w) = \frac{-(1+\sqrt{5}) \eta^{\tau_5} + 2}{2\eta^{\tau_5}+1+\sqrt{5}},$$
(for all $-d=d_Kf^2 <-4$) where $\tau_5 = \left(\frac{\mathbb{Q}(\eta)/K}{\wp_5}\right)$ is the Frobenius automorphism (Artin symbol) for $\wp_5$ (which is defined since $\mathbb{Q}(r(w/5))$ is abelian over $K$ and unramified at $\wp_5$).
See Tables 1 and 2 for a list of the polynomials $p_d(x)$ for small values of $d$.  As is clear from the tables, these polynomials have relatively small coefficients and discriminants.  Moreover, as we show in Section 5, these values of $r(\tau)$ are periodic points of an algebraic function, and can be computed for small values of $d$ and small periods using nested resultants.  (See  \cite{mor5}, Section 3, and  \cite{mor7}.)  We prove the following.  \bigskip

\noindent {\bf Theorem 1.2.} {\it If
$$g(X,Y)=(Y^4+2Y^3+4Y^2+3Y+1)X^5-Y(Y^4-3Y^3+4Y^2-2Y+1),$$
the roots of $p_d(x)$ are periodic points of the multi-valued algebraic function $\mathfrak{g}(z)$ defined by $g(z,\mathfrak{g}(z))=0$.  With $w$ chosen as in Theorem 1.1, the period of $\eta=r(w/5)$ with respect to the action of $\mathfrak{g}$ is the order of the Frobenius automorphism $\displaystyle \tau_5=\left(\frac{\mathbb{Q}(\eta)/K}{\wp_5} \right)$ in $\textrm{Gal}(\mathbb{Q}(\eta)/K)$.} \bigskip

As part of our discussion we also prove the following. To state the result, let
$$\mathfrak{s}(z) = \frac{(\zeta+\zeta^2)z+1}{z+1+\zeta+\zeta^2}, \ \ \zeta=\zeta_5=e^{2\pi i/5},$$
a linear fractional map of order $5$.  The group $\langle \mathfrak{s}(z) \rangle$ generated by $\mathfrak{s}(z)$ is the Galois group of the function field extension $\mathbb{Q}(\zeta,z)/\mathbb{Q}(\zeta, \mathfrak{r}(z))$, where
$$\mathfrak{r}(z)=\frac{z(z^4-3z^3+4z^2-2z+1)}{z^4+2z^3+4z^2+3z+1}.$$

\noindent {\bf Theorem 1.3.}  {\it With $w$ as in Theorem 1.1 and $\tau_5$ as above, we have the formula
$$r(w/5)^{\tau_5} = \mathfrak{s}^j(r(w))=r\left(\frac{w}{1-jw}\right),$$
where $j \not \equiv 0$ (mod $5$) has the same value as in Theorem 1.1 and $j$ is the unique integer (mod $5$) for which $\mathfrak{s}^j(r(w))$ is an algebraic conjugate of $\eta=r(w/5)$.} \bigskip

This fact is significant, because in the ideal-theoretic formulations of Shimura's Reciprocity Law, such as in \cite{sch}, p. 123, one has to restrict to ideals that are relatively prime to the level of the modular function being considered.  Here $r(\tau) \in \Gamma(5)$, so the level is $N=5$, but Theorem 1.3 gives information about the automorphism $\tau_5$ corresponding to the prime ideal $\wp_5$ of $K$. \medskip

Theorem 1.3 has the following application.  A formula for the real value
$$r(3i)=\frac{e^{-6 \pi/5}}{1+} \ \frac{e^{-6 \pi}}{1+} \ \frac{e^{-12 \pi}}{1+} \ \frac{e^{-18 \pi}}{1+} \dots$$
was stated by Ramanujan in his notebooks and proved in \cite{bch} and \cite{bcz}.  In Section 5 we prove the alternative formula
$$r(3i)=\frac{(1+\zeta^3)\eta^{\tau_5}+\zeta}{\eta^{\tau_5}-\zeta-\zeta^3},\eqno{(1.6)}$$
where
$$\eta^{\tau_5}=r\left(\frac{4+3i}{5}\right)^{\tau_5} =\frac{-i \omega}{2}-\frac{i\sqrt{3}}{2}+i\frac{\omega^2}{4} \sqrt[4]{3}\left(\sqrt{4+2\sqrt{5}}+i\sqrt{-4+2\sqrt{5}}\right)$$
and $\omega=(-1+i \sqrt{3})/2$.  This formula expresses Ramanujan's value in terms of roots of unity and simpler square-roots than appear in his original formula.  (See Example 1 in Section 5.)  \medskip

\section{Defining the Heegner points.}

Throughout the paper we will have occasion to make use of the linear fractional map
$$\tau(b)=\frac{-b+\varepsilon^5}{\varepsilon^5 b+1}=\frac{-b+\varepsilon_1}{\varepsilon_1 b+1}, \quad \varepsilon_1=\varepsilon^5=\frac{-11+5\sqrt{5}}{2}. \eqno{(2.1)}$$
We note that
$$j_5(\tau(b))=j_{5,5}(b)=\frac{(1+228b+494b^2-228b^3+b^4)^3}{b(1-11b-b^2)^5},\eqno{(2.2a)}$$
$$=-\frac{(z^2-228z+496)^3}{(z+11)^5}, \ \ z=b-\frac{1}{b},\eqno{(2.2b)}$$
where $j_{5,5}(b)$ is the $j$-invariant of the elliptic curve
\begin{align*}
E_{5,5}(b): \ Y^2+(1+b)XY+5bY=&X^3+7bX^2+(6b^3+6b^2-6b)X\\
&+ \ b^5+b^4-10b^3-29b^2-b.
\end{align*}
The curve $E_{5,5}(b)$ is isogenous to $E_5(b)$ (\cite{mor}, p. 259), and because of (2.2), $E_5(\tau(b))$ represents the Tate normal form for $E_{5,5}(b)$. \medskip

Let $K=\mathbb{Q}(\sqrt{-d})$, where $-d=d_Kf^2 \equiv \pm 1$ (mod $5$) and $d_K$ is the discriminant of $K$.  As usual, let $\eta(\tau)$ be the Dedekind $\eta$-function.  From Weber  \cite{w}, p.256, the function
$$x_1=x_1(w)=\left(\frac{\eta(w/5)}{\eta(w)}\right)^2$$
satisfies the equation
$$x_1^6+10x_1^3-\gamma_2(w)x_1+5=0, \ \ \gamma_2(w)=j(w)^{1/3}.$$
Thus
$$j(w)=\frac{(x_1^6+10x_1^3+5)^3}{x_1^3}.\eqno{(2.3)}$$
On the other hand,
$$x_1^3=y^5+5y^4+15y^3+25y^2+25y=(y+1)^5+5(y+1)^3+5(y+1)-11,$$
with $y=y(w)=\frac{\eta(w/25)}{\eta(w)}$.  By Theorem 6.6.4 of Schertz  \cite{sch}, p. 159, both $x_1^3$ and $y$ are elements of the ring class field $\Omega_f=K(j(w))$ if 
$$w=\begin{cases}
\frac{v+\sqrt{-d}}{2}, &2 \nmid d, \ v^2 \equiv -d \ (\textrm{mod} \ 5^2), \ (v,2f)=1,\\
v + \frac{\sqrt{-d}}{2}, &2 \mid d, \ 2 \nmid f, \ v^2 \equiv -d/4 \ (\textrm{mod} \  5^2), \ \ (v,f)=1,\\
v + \frac{\sqrt{-d}}{2}, &2 \mid d, \ 2 \mid f, \ v^2 \equiv -d/4 \ (\textrm{mod} \  5^2), \ \ (v,f_{odd})=1;
\end{cases} \eqno{(2.4)}$$
in the last case $f_{odd}$ is the largest odd divisor of $f$ and $v \not \equiv d/4$ (mod $2$) is chosen to guarantee that $(N(w), f)=1$.  (The latter condition is needed to insure that $(w)$ is a proper ideal of $\textsf{R}_{-d}$ in Section 4.)  These conditions on $w$ are equivalent to the conditions imposed on $w$ in Theorem 1.1.  \medskip

Now we set
$$z=z(w)=b-\frac{1}{b}=-11-x_1^3=-11-\left(\frac{\eta(w/5)}{\eta(w)}\right)^6,\eqno{(2.5)}$$
so that $b$ is one of the two roots of the equation
$$b^2-zb-1=0, \ \ z=-11-x_1^3.$$
From the identity
$$\frac{1}{r^5(\tau)}-11-r^5(\tau)=\left(\frac{\eta(\tau)}{\eta(5\tau)}\right)^6, \ \ \tau \in \mathbb{H},$$
for the Rogers-Ramanujan function $r(\tau)$ (see  \cite{du}), we see that
$$\frac{1}{b}-b -11 = \frac{1}{r^5(w/5)}-r^5(w/5)-11,$$
from which it follows that
$$b=r^5(w/5) \ \ \textrm{or} \ \ \frac{-1}{r^5(w/5)}\eqno{(2.6)}$$
and
$$z=r^5(w/5)-\frac{1}{r^5(w/5)}.\eqno{(2.7)}$$
We find from (1.5), (2.5), and (2.3) that
$$j_5(b)=\frac{((-11-x_1^3)^2+12(-11-x_1^3)+16)^3}{x_1^3}$$
$$=\frac{(x_1^6+10x_1^3+5)^3}{x_1^3}=j(w).\eqno{(2.8)}$$
When $z$ is given by (2.5), $j(w)$ is the $j$-invariant of $E_5(b)$.  Weber  \cite{w}, p.256, also gives the equation
$$j(w/5)=\frac{(x_1^6+250x_1^3+3125)^3}{x_1^{15}}=j_{5,5}(b),\eqno{(2.9)}$$
for the same substitution (2.5), by (2.2{\it b}).  Thus, $j(w/5)$ is the $j$-invariant of the isogenous curve $E_{5,5}(b)$. \medskip

The functions $z(w)$ and $y(w)$ are modular functions for the group $\Gamma_0(5)$, by Schertz  \cite{sch}, p. 51.  Moreover, $w$ and $w/5$ are basis quotients for proper ideals in the order $\textsf{R}_{-d}$ of discriminant $-d$ in $K$.  Hence, we have the following.  \smallskip

\noindent {\bf Theorem 2.1.}  {\it If $z=b-1/b$ satisfies (2.5), where $w$ is given by (2.4), then $j_5(b)=j(w)$ and $j_{5,5}(b)=j(w/5)$ are roots of the class equation $H_{-d}(x)=0$, and the isogeny $E_5(b) \rightarrow E_{5,5}(b)$ represents a Heegner point on $\Gamma_0(5)$.  Furthermore, $z$ lies in the ring class field of conductor $f$ over $K=\mathbb{Q}(\sqrt{-d})$, where $-d=f^2 d_K$ and $d_K$ is the discriminant of $K$.}  \bigskip

Exactly the same arguments apply if $w$ is replaced in (2.3)-(2.9) by $w/a$, where $(a,f)=1$ and $5a \mid N(w)$.  (To guarantee $y(w/a) \in \Omega_f$ we would also need $5^2 a \mid N(w)$.)  Then $w/a$ and $w/(5a)$ are basis quotients for proper ideals in $\textsf{R}_{-d}$ and $j(w/a)$ and $j(w/(5a))$ are roots of $H_{-d}(x)$.  Thus, $j(w), j(w/a) \in \Omega_f$ are conjugate to each other over $K$.  Theorem 6.6.4 of Schertz \cite{sch} implies that the corresponding values $z(w), z(w/a)$ in (2.5) are also conjugate to each other over $K$ if $5 \nmid a$, but in Section 4 we will need to relax this restriction on $a$.  To do this, we prove the following lemma.  Let $J(z)$ denote the rational function
$$J(z)=-\frac{(z^2+12z+16)^3}{z+11}.$$
Recall that an ideal $\mathfrak{a}$ of the order $\textsf{R}_{-d}$ corresponds to the ideal $\mathfrak{a} R_K$ of the maximal order $R_K = \textsf{R}_{d_K}$ of $K$, and conversely, an ideal $\mathfrak{b}$ in $R_K$ corresponds to the ideal $\mathfrak{b}_d=\mathfrak{b} \cap \textsf{R}_{-d}$ in $\textsf{R}_{-d}$ (\cite{co}, p. 144).  \bigskip

\noindent {\bf Lemma 2.2.} {\it For a given ideal $\mathfrak{a} = (a,w) \subseteq \textsf{R}_{-d}$ with ideal basis quotient $w/a$, where $(a,f)=1$ and $5a \mid N(w)$, there is a unique value of $z_1 \in \Omega_f$ for which $J(z_1) = j(w/a)$ and $z_1+11 \cong \wp_5'^3$, and this value is $z_1 = z^{\sigma^{-1}}$, where $\sigma = \left(\frac{\Omega_f/K}{\mathfrak{a}R_K}\right)$.  ($\alpha \cong$ $\beta$ denotes equality of the divisors $(\alpha)$ and $(\beta)$.)}  \medskip

\noindent {\it Proof.}  If $\sigma$ is the Frobenius automorphism given in the statement of the lemma, $j(w/a)^\sigma=j(\mathfrak{a})^\sigma = j(\textsf{R}_{-d}) = j(w) = J(z)$, it follows that $J(z^{\sigma^{-1}}) = j(w/a)$.  Suppose there is a $z_2 \in \Omega_f$, different from $z_1 = z^{\sigma^{-1}}$, for which $J(z_2)=J(z_1)$ and $z_2+11 \cong z_1+11$.  Then $(z_1,z_2)$ is a point on the curve $F(u,v)=0$, where
\begin{align*}
F(u,v) &= -(u+11)(v+11)\frac{J(u)-J(v)}{u-v}\\
&=(v+11)u^5+(v^2+47v+396)u^4+(v^3+47v^2+876v+5280)u^3\\
& \ +(v^4+47v^3+876v^2+8160v+31680)u^2\\
& \ +(v^5+47v^4+876v^3+8160v^2+39360v+84480)u\\
& \ +11v^5+396v^4+5280v^3+31680v^2+84480v+97280.
\end{align*}
A calculation on Maple shows that this is a curve of genus $0$, parametrized by the rational functions
\begin{align*}
u&=-\frac{11t^5+55t^4+165t^3+275t^2+275t+125}{t(t^4+5t^3+15t^2+25t+25)}\\
v&=-\frac{t^5+11t^4+55t^3+165t^2+275t+275}{t^4+5t^3+15t^2+25t+25}.
\end{align*}
Hence, $F(z_1,z_2)=0$ gives that
$$z_1+11= \frac{-125}{t(t^4+5t^3+15t^2+25t+25)},$$
or
$$t^5+5t^4+15t^3+25t^2+25t+\frac{125}{z_1+11}=0,$$
for some algebraic number $t$.  Since $z_1+11 \cong z+11 \cong \wp_5'^3$ (see eq. (4.2) below), we have $(z_1+11) \mid 5^3$ and $t$ is an algebraic integer which is not divisible by any prime divisor of $\wp_5'$ in $\Omega_f(t)$.  Then
$$z_2+11=\frac{-t^5}{t^4+5t^3+15t^2+25t+25}=\frac{t^5}{\frac{125}{t(z_1+11)}}=t^6 \frac{(z_1+11)}{125}.$$
 But the equality of the ideals $(z_2+11) = (z_1 +11)$ implies that $t^6 \cong 5^3$, so $t$ {\it is divisible} by some prime divisor of $\wp_5'$ in $\Omega_f(t)$.  This contradiction establishes the claim.  $\square$

\section{Points of order 5 on $E_5(b)$.}

From  \cite{mor6} we take the following.  The $X$-coordinates of points of order 5 on $E_5(b)$ which are not in the group 
$$\langle (0,0) \rangle = \{O, (0,0), (0,-b), (-b,0), (-b,b^2)\}$$
can be given in the form
\begin{align*}
X&=\frac{(5-\alpha)}{100} \{(-18-12b+6b\alpha+8\alpha-2b^2)u^4+(-4b\alpha+2b^2+3\alpha-7+12b)u^3\\
&+(7b\alpha+\alpha-3-2b^2-7b)u^2+(22b-2+2b^2)u-3-7b+3b\alpha-2b^2-\alpha \}\\
&=\frac{(5-\alpha)}{100}(A_4u^4+A_3u^3+A_2u^2+A_1u+A_0),
\end{align*}
where $\alpha=\pm \sqrt{5}$, 
$$u^5=\phi_1(b) = \frac{2b+11+5\alpha}{-2b-11+5\alpha}=\frac{b-\bar \varepsilon^5}{-b+\varepsilon^5} \eqno{(3.1)}$$
and
$$\varepsilon = \frac{-1+\alpha}{2}, \ \ \bar \varepsilon = \frac{-1-\alpha}{2}.$$
Equation (3.1) shows that $u^5 = 1/(\varepsilon^5 \tau(b))$, i.e., $\tau(b) = (\varepsilon u)^{-5}$.  Solving for $b$ in (3.1) gives
$$b = \frac{\varepsilon^5u^5+\bar \varepsilon^5}{u^5+1}.\eqno{(3.2)}$$ \smallskip

\noindent Now the Weierstrass normal form of $E_5(b)$ is given by
$$Y^2=4X^3-g_2X-g_3,\ \ g_2=\frac{1}{12}(b^4+12b^3+14b^2-12b+1),$$
$$g_3=\frac{-1}{216}(b^2+1)(b^4+18b^3+74b^2-18b+1),$$
with
$$\Delta=g_2^3-27g_3^2=b^5(1-11b-b^2).$$ \medskip
By Theorem 2.1, $E_5(b)$ has complex multiplication by the order $\textsf{R}_{-d}$, so the theory of complex multiplication implies that if $K \neq \mathbb{Q}(i)$, i.e. $d \neq 4f^2$, the $X$-coordinates $X(P)$ of points of order $5$ on $E_5(b)$ have the property that the quantities
$$\frac{g_2g_3}{\Delta}\left(X(P)+\frac{1}{12}(b^2+6b+1)\right)$$
generate the field $\Sigma_5 \Omega_f$ over $\Omega_f$, where $\Sigma_5$ is the ray class field of conductor $5$ over $K=\mathbb{Q}(\sqrt{-d})$.  (See  \cite{fra}; or  \cite{si} for $f=1$.)  \medskip

In the case that $d=4f^2>4$, the argument leading to Theorem 2 of  \cite{fra} shows that these quantities generate a class field $\Sigma'_{5f}$ over $K=\mathbb{Q}(i)$ whose corresponding ideal group $\textsf{H}$ consists of the principal ideals generated by elements of $K$, prime to $5f$, which are congruent to rational numbers (mod $f$) and congruent to $\pm 1$ (mod 5). $\textsf{H}$ is an ideal group because it contains the ray mod $5f$.  Thus $\textsf{H} \subset \textsf{S}_5 \cap \textsf{P}_f$ is contained in the intersection of the principal ring class mod $f$, $\textsf{P}_f$, and the ray mod $5$, $\textsf{S}_5$.  If $(\alpha) \in \textsf{S}_5 \cap \textsf{P}_f$, then we may take $\alpha \equiv r$ (mod $f$) and $r \in \mathbb{Q}$; and then $i^a \alpha \equiv 1$ (mod $5$) for some power of $i$.  If $2 \mid a$, then $(\alpha) \in \textsf{H}$; while if $2 \nmid a$, then $\alpha^2 \equiv -1$ (mod $5$), so $(\alpha)^2 \in \textsf{H}$, and the product of any two such ideals lies in $\textsf{H}$.  This implies that $[\textsf{S}_5 \cap \textsf{P}_f: \textsf{H}]=2$ and $\Sigma'_{5f}$ is a quadratic extension of $\Sigma_5\Omega_f$ (when $K=\mathbb{Q}(i)$).  Moreover, $\textsf{H}$ is a subgroup of the principal ring class $\textsf{P}_{5f}$ and $[\textsf{P}_{5f}:\textsf{H}]=2$, so that $[\Sigma_{5f}' : \Omega_{5f}]=2$.  Since $\textsf{P}_{5f} \neq \textsf{S}_5 \cap \textsf{P}_f$, it follows that $\Sigma_{5f}' = \Omega_{5f} (\Sigma_5 \Omega_f) = \Sigma_5\Omega_{5f}$.  Noting that $\textsf{P}_f/\textsf{P}_{5f}$ is cyclic of order $4$, generated by $(\alpha)\textsf{P}_{5f}$ with $\alpha \equiv 2$ (mod $\wp_5$) and $\equiv 1$ (mod $\wp_5'$), it follows from Artin Reciprocity that $\Omega_{5f}/\Omega_f$ is a cyclic quartic extension.  \medskip

Let $F$ denote the field $\Sigma_5\Omega_f$, for $d \neq 4f^2$; and $\Sigma_{5f}'=\Sigma_5\Omega_{5f}$, for $d=4f^2>4$.  Also, let $\phi(\mathfrak{a})$ denote the Euler $\phi$-function for ideals $\mathfrak{a}$ of $R_K$.  Since $p=5=\wp_5 \wp_5'$ splits in $K$, the degree of $\Sigma_5/\Sigma_1$ is given by
$$[\Sigma_5: \Sigma_1]=\frac{1}{2}\phi(\wp_5)\phi(\wp_5')=8, \ \ \textrm{if} \ d \neq 4f^2;$$
and since every intermediate field of $\Sigma_5/\Sigma_1$ is ramified over $p=5$ we have that
$$[F:\Omega_f]=[\Sigma_5\Omega_f: \Omega_f]=8, \ \ d \neq 4f^2.$$
On the other hand,
$$[F:\Omega_f]=[\Sigma'_{5f}:\Omega_f]=2 \cdot [\Sigma_5\Omega_f:\Omega_f]=8, \ \ d=4f^2>4,$$
since in this case
$$[\Sigma_5: K]=\frac{1}{4}\phi(\wp_5)\phi(\wp_5')=4, \ \ d=4f^2;$$
so that $\Sigma_5=K(\zeta_5)$ when $K=\mathbb{Q}(i)$.  Thus, $[F:\Omega_f]=8$ in all cases (with $d \neq 4$). \medskip

We henceforth take $\alpha=\sqrt{5}$ in the above formulas, and we prove the following. \bigskip

\noindent {\bf Theorem 3.1.}  {\it If $z=b-1/b$ is given by (2.7), where $w$ is given by (2.4), with $d \neq 4$, then the roots $u$ of the equation (3.1) lie in the field $F=\Sigma_5 \Omega_f$, if $d \neq 4f^2$, and in $F=\Sigma_5\Omega_{5f}$, if $d=4f^2>4$. Thus, the value $b$ is given by
$$b=\frac{\varepsilon^5 u^5+ \bar \varepsilon^5}{u^5+1}, \ \ \varepsilon = \frac{-1+\sqrt{5}}{2}, \ \ \bar \varepsilon=\frac{-1-\sqrt{5}}{2},$$
where
$$u = -\frac{r(w)-\bar \varepsilon}{r(w)-\varepsilon} \ \ \textrm{or} \ \ -\frac{\bar \varepsilon r(w)+1}{\varepsilon r(w)+1},$$
according as $b=r^5(w/5)$ or $b=\frac{-1}{r^5(w/5)}$.  Moreover,  $r(w), r(w/5)$ and $r(-1/w)$ lie in the field $F$.}  \medskip 

\noindent {\it Proof.}  Note first that
$$\frac{g_2g_3}{\Delta}=\frac{-1}{2592}\frac{(b^4+12b^3+14b^2-12b+1)(b^2+1)(b^4+18b^3+74b^2-18b+1)}{b^5(1-11b-b^2)}$$
$$=\frac{1}{2592} \frac{(z^2+12z+16)(z^2+18z+76)}{z+11} \frac{b^2+1}{b^2},$$ \smallskip

\noindent where $z=b-\frac{1}{b}=-11-x_1^3$ lies in $\Omega_f$.  It follows that
$$\frac{b^2+1}{b^2}\left(X(P) + \frac{1}{12}(b^2+6b+1)\right) \in F$$
for any point $P \in E_5[5]$.  In particular, with $P=(-b,0)$ we have that
$$\frac{b^2+1}{12b^2}(b^2-6b+1) = \frac{1}{12}\left(b+\frac{1}{b}\right)\left(b+\frac{1}{b}-6 \right) \in F.$$
Since $b-\frac{1}{b}$ lies in $\Omega_f$, the field $F$ contains the quantity
$$\left(b-\frac{1}{b} \right)^2+4=b^2+\frac{1}{b^2}+2=\left(b+\frac{1}{b} \right)^2,$$
and therefore also $\left(b+\frac{1}{b}\right)$ and $\left(b+\frac{1}{b}\right) + \left(b-\frac{1}{b}\right) = 2b$.  Therefore, $b \in F$ and we have that
$$X(P) \in F, \ \ \textrm{for} \ P \in E_5[5].$$

Since $\mathbb{Q}(\sqrt{5}) \subset \mathbb{Q}(\zeta_5) \subseteq \Sigma_5$, we deduce from the formula for $X$ above that
$$A_4u^4+A_3u^3+A_2u^2+A_1u+A_0 \in F$$
for any root of (3.1).  Hence, for any fixed root $u$ of (3.1) we have that
$$A_4\zeta^{4i}u^4+A_3\zeta^{3i}u^3+A_2\zeta^{2i}u^2+A_1\zeta^{i}u+A_0 = B_i \in F, \ \ 0 \le i \le 4.\eqno{(3.3)}$$
This gives a system of 5 equations in the 5 ``unknowns'' $u^i$, with coefficients in $F$.  The determinant of this system is
$$D=-\frac{5^2}{8}(\zeta-\zeta^2-\zeta^3+\zeta^4)(-3-7b+3b\alpha-2b^2-\alpha)(-2b-1+\alpha)(2b+\alpha+1)$$
$$\times (2b+11+5\alpha)(-b+2+\alpha)(-2b-11+5\alpha)^4,\eqno{(3.4)}$$
which I claim is not zero. \medskip

Ignoring the constant term $\frac{\pm 5^2 \sqrt{5}}{8}$ in front, multiply the rest by the polynomial in (3.4) obtained by replacing $\alpha$ with $-\alpha$.  This gives the polynomial
$$2^{16}(b^2-4b-1)(b^4+7b^3+4b^2+18b+1)(b^2+11b-1)^5(b^2+b-1)^2.$$
If $b$ is a root of any of the quadratic factors, then $z=b-\frac{1}{b}$ is rational: $z=4, -11$, or $-1$, respectively.  In these cases $j(w)=-102400/3, \infty$, or $-25/2$, all of which are impossible, since $j(w)$ is an algebraic integer.  \medskip

Now $E_5(b)$ has complex multiplication by an order in the field $K=\mathbb{Q}(\sqrt{-d})$ whose discriminant is not divisible by $5$.  Therefore, $j(w)=j(E_5(b))$ generates an extension of $\mathbb{Q}$ which is not ramified at $p=5$.  If $b$ is a root of $h(x)=x^4+7x^3+4x^2+18x+1$, then $\textrm{disc}(h(x)) = -5^8 19$ and $\textrm{Gal}(h(x)/\mathbb{Q}) \cong D_4$ imply that $K(b)$ can only be abelian over the quadratic field $K=\mathbb{Q}(\sqrt{-19})$ and $f=1$.  Then $j_5(b)$ is a root of the irreducible polynomial
$$H(x)=x^4+5584305x^3-32305549025x^2+63531273863125x-5^6 31^3 449^3,$$
which is impossible, since $K=\mathbb{Q}(\sqrt{-19})$ has class number $1$.  This shows that the determinant $D$ in (3.4) is nonzero, and therefore, since the coefficients $A_i$ and $D$ lie in the field $F$, we get that the solution $(u^4,u^3,u^2,u,1)$ of the system (3.3) lies in $F$ also.  This proves that $u \in F$.  In particular, $\tau(b) = (\varepsilon u)^{-5}$ is a $5$-th power in $F$.  \medskip

We can find formulas for $u$ from the identities
$$r^5\left(\frac{-1}{5\tau}\right) =\frac{-r^5(\tau)+\varepsilon^5}{\varepsilon^5 r^5(\tau)+1} \ \ \textrm{and} \ \ r\left(\frac{-1}{5\tau}\right)= \frac{\bar \varepsilon r(5\tau)+1}{r(5\tau)-\bar \varepsilon}.$$
See  \cite{du}, p. 150.  If $\tau=w/5$ and $b=r^5(w/5)$, we have
$$r^5\left(\frac{-1}{w}\right) = \frac{-b+\varepsilon^5}{\varepsilon^5 (b-\bar \varepsilon^5)}=\frac{1}{\varepsilon^5 u^5},$$
and we can take
$$u = \frac{1}{\varepsilon r\left(\frac{-1}{w}\right)}= \frac{r(w)-\bar \varepsilon}{\varepsilon(\bar \varepsilon r(w)+1)}=-\frac{r(w)-\bar \varepsilon}{r(w)-\varepsilon}, \ \ b=r^5(w/5).\eqno{(3.5)}$$
On the other hand, if $b=\frac{-1}{r^5(w/5)}$, then we can choose
$$u = -\frac{\bar \varepsilon r(w)+1}{\varepsilon r(w)+1}.\eqno{(3.6)}$$
In either case it is clear that $r(w), r(-1/w) \in F$.  \medskip

We can apply the same analysis with $b$ replaced by $\tau(b)$, since $E_{5,5}(b) \cong E_5(\tau(b))$, so that the latter curve also has complex multiplication by $\textsf{R}_{-d}$.  Furthermore,
$$b=r^5(w/5) \ \Longrightarrow \ \tau(b) = r^5\left(\frac{-1}{w}\right),$$
while
$$b=\frac{-1}{r^5(w/5)} \ \Longrightarrow \ \tau(b) = \frac{-1}{r^5(-1/w)}.$$
Note also that when $b$ is replaced by $\tau(b)$ in the determinant $D$, its factors in $b$ are
$$\frac{(2b+1)(b-2)(b+3)(-3-7b+3b\alpha-2b^2-\alpha)b^4}{(2b+11+5\alpha)^{10}},$$
and so are nonzero by the same reason as before.  Hence we get a solution $u_1 \in F$ of the equation
$$u_1^5=\phi_1(\tau(b))=-\frac{\bar \varepsilon^5}{b}=\frac{1}{\varepsilon^5 b}.$$
Therefore, $b=1/(\varepsilon u_1)^5$ is also a 5-th power in $F$, i.e. $r(w/5) \in F$. $\square$  \bigskip

\noindent {\bf Remarks.} (1) The result of Theorem 3.1 that $r(w), r(w/5) \in F$ is sharper than what is obtained from  \cite{sch}, Thm. 5.1.2, p. 123.  That theorem only yields that $r(w), r(w/5) \in \Sigma_{5f}$, the ray class field of conductor $5f$.  Also, the coefficients of the $q$-expansion of $r(-1/\tau)$ are in $\mathbb{Q}(\sqrt{5})$ but not all in $\mathbb{Q}$, so \cite{sch}, Theorem 5.2.1 does not apply.  \smallskip

\noindent (2) The results of  \cite{mor6} show that the coordinates of all the points in $E_5(b)[5]-\langle (0,0) \rangle$ are rational functions of the quantity $u$, and therefore of the quantity $r(w)$, with coefficients in $\mathbb{Q}(\zeta_5)$, by (3.5) and (3.6).  It follows from the theory of complex multiplication that $F=K(\zeta_5, r(w))$.  In Corollary 4.7 below we will prove that $F=\mathbb{Q}(r(w))$.  \bigskip

Now $b$ satisfies the equation $b-\frac{1}{b}=z=-11-x_1^3 \in \Omega_f$, so $b$ is at most quadratic over $\Omega_f$.  Hence, its degree over $\mathbb{Q}$ is at most $4h(-d)$.  This degree is also at least $h(-d)$ since $j(w) \in \mathbb{Q}(b)$.  \bigskip

\noindent {\bf Proposition 3.2.} {\it If $d > 4$, the degree of $z=b-1/b$ over $\mathbb{Q}$ is $2h(-d)$.  Thus, $\Omega_f=\mathbb{Q}(z)$, and the minimal polynomial $R_d(X)$ of $z$ over $\mathbb{Q}$ is normal.} \smallskip

\noindent {\it Proof.}  Recall from above that
$$j(w)=j_5(b)=-\frac{(z^2+12z+16)^3}{z+11},$$
and
$$j(w/5)=j_{5,5}(b)=-\frac{(z^2-228z+496)^3}{(z+11)^5}.$$
Since $z=-11-x_1^3 \in \Omega_f$ and the real number $j(w)$ has degree $h(-d)$ over $\mathbb{Q}$, it is clear that the degree of $z$ is either $h(-d)$ or $2h(-d)$.  Suppose the degree is $h(-d)$.  Then $\mathbb{Q}(z)=\mathbb{Q}(j(w))$, which implies that $z$ is real, and therefore $j(w/5)$ is also real.  We also know $j(w/5)=j(\wp_{5,d})$, where $\wp_{5,d} =\wp_5 \cap \textsf{R}_{-d}$, so that $j(\wp_{5,d}) = \overline{j(\wp_{5,d})}=j(\wp_{5,d}^{-1})$ implies that $\wp_5$ must have order $1$ or $2$ in the ring class group of $K$ (mod $f$). \medskip

If $\wp_5 \sim 1$ (mod $f$), then $4\cdot 5=x_2^2+dy_2^2$ for some integers $x_2, y_2$, which implies that $d = 4, 11, 16, 19$, the first of which is excluded.  In the last three cases we have, respectively
$$H_{-11}(x)=x+32^3, \ \ H_{-16}(x)=x-66^3, \ \ H_{-19}(x)=x+96^3.$$
(See  \cite{co}.)  In these cases there is only one irreducible polynomial $Q_d(x)$ of degree $4h(-d)=4$ or less which divides $F_d(x)$ in (1.2), which must therefore be the minimal polynomial of $b$.  We have
$$Q_{11}(x)=x^4+4x^3+46x^2-4x+1, \ \ Q_{16}(x) = x^4+18x^3+200x^2-18x+1,$$
$$Q_{19}(x) = x^4+36x^3+398x^2-36x+1.$$
To each of these polynomials with root $b$ corresponds the minimal polynomial $R_d(x)$ with root $z=b-\frac{1}{b}$.  These are:
$$R_{11}(x)=x^2+4x+48, \ \  R_{16}(x) = x^2+18x+202, \ \ R_{19}(x)=x^2+36x+400,$$
each of which has the correct degree $2h(-d)=2$. \medskip

Now suppose that the order of $\wp_5$ is $2$.  Then $\wp_5^2 \sim 1$ (mod $f$) implies that $4\cdot 5^2=x_2^2+dy_2^2$ for $x_2, y_2 \in \mathbb{Z}$ with $x_2 \equiv y_2$ (mod $2$), if $d$ is odd, giving the possibilities:
$$d=51, 91, 99, \ \ \textrm{with} \ \ h(-51)=h(-91)=h(-99)=2;$$
and $5^2=x_2^2+\frac{d}{4}y_2^2$, if $d$ is even, in which case we have the following possibilities:
$d = 24, 36, 64, 84, 96$, with
$$h(-24)=h(-36)=h(-64)=2, \ \ h(-84)=h(-96)=4.$$
We use the following class equations (see Fricke \cite{fr1}, III, pp. 401, 405, 420 for $D=-24, -36, -64, -91$; and Fricke  \cite{fr2}, III, p. 201 for $D=-51$):
$$H_{-24}(x)=x^2-4834944x+14670139392,$$
$$H_{-36}(x)=x^2-153542016x-1790957481984,$$
$$H_{-51}(x)=x^2+5541101568x+6262062317568,$$
$$H_{-64}(x)=x^2-82226316240x-7367066619912,$$
$$H_{-91}(x)=x^2+10359073013760x-3845689020776448,$$
$$H_{-99}(x)=x^2+37616060956672x-56171326053810176.$$
These polynomials yield the following minimal polynomials for $z$:
$$R_{24}(x)=x^4-12x^3+20x^2+3120x+16912,$$
$$R_{36}(x)=x^4+60x^3+3020x^2+51984x+287248,$$
$$R_{51}(x)=x^4-24x^3+6800x^2+155136x+852736,$$
$$R_{64}(x)=x^4-216x^3+17234x^2+430380x+2362354,$$
$$R_{91}(x)=x^4-216x^3+154448x^2+3449088x+18965248,$$
$$R_{99}(x)=x^4+872x^3+292624x^2+6230016x+34284288.$$
We computed $H_{-99}(x)$ and $R_{99}(x)$ directly from (2.5).  In the same way we find
$$R_{84}(x)=x^8-468x^7+81124x^6+3053232x^5+65642496x^4+1156633920x^3$$
$$+13586087488x^2+88268813568x+244368064768,$$ \smallskip
$$R_{96}(x)=x^8+324x^7+230848x^6+5080248x^5+32351604x^4+88662672x^3$$
$$+675333328x^2+2681910144x+7697193232.$$
Each of these polynomials is irreducible, so the quantity $z$ always has degree $2h(-d)$ over $\mathbb{Q}$.  Since $z \in \Omega_f$, it follows that $\Omega_f=\mathbb{Q}(z)$.  This proves the claim.  $\square$  \bigskip

\noindent {\bf Theorem 3.3.}  {\it With $z$ as in (2.7) and $d>4$, the quantities $b$ and $\displaystyle \tau(b) = \frac{-b+\varepsilon^5}{\varepsilon^5 b+1}$
are $5$-th powers in the field $F$, and if
$$\xi^5=\tau(b) \ \ \textrm{and} \ \ \eta^5=b,\eqno{(3.7)}$$
then $(X,Y)=(\xi, \eta)$ is a solution in $F$ of the equation
$$X^5+Y^5=\varepsilon^5(1-X^5Y^5).\eqno{(3.8)}$$
Such numbers $\xi$ and $\eta$ exist for which $\xi \in \mathbb{Q}(\tau(b))$ and $\eta \in \mathbb{Q}(b)$.}  \bigskip

\noindent {\it Proof.} From (3.7) and the last part of the proof of Proposition 3.1, we have
$$b=\frac{1}{\varepsilon^5 u_1^5}=\eta^5, \ \ \tau(b) = \frac{1}{\varepsilon^5 u^5}=\xi^5;$$
with
$$\eta=\delta \zeta^i r^{\delta}\left(\frac{w}{5}\right), \ \ \xi = \delta \zeta^{\delta j} r^{\delta}\left(\frac{-1}{w}\right), \ \ \delta=\pm 1.\eqno{(3.9)}$$
The relation $\xi^5 = \tau(\eta^5)$ implies that $(X,Y)=(\xi,\eta)$ lies on (3.8).  It only remains to prove that $\eta=\frac{1}{\varepsilon u_1}=b^{1/5}$ can be chosen to lie in $\mathbb{Q}(b)$.  The polynomial $q(X)=X^5-b$ has the root $\eta$ and splits completely in $F$.  Since the degree $[F: \Omega_f]=8$ is not divisible by $5$ or by $3$, and the degree $[\mathbb{Q}(b):\Omega_f]=[\mathbb{Q}(b):\mathbb{Q}(z)]$ divides $2$, $q(X)$ has to factor into a product of a linear and a quartic polynomial, or a linear times a product of two quadratics over $\mathbb{Q}(b)$.  Hence, at least one root of $q(X)$ has to lie in $\mathbb{Q}(b)$, and we can assume this root is $\eta$.  In the same way, we can assume $\xi \in \mathbb{Q}(\tau(b))$.  $\square$  \bigskip

\noindent {\bf Remark.}  When $d=4$, $(X,Y)=(\xi,\eta)=(-i,i)$ is a solution of the equation (3.8), corresonding to the values $b=i, z=2i$. \medskip

Using (3.7), we see that
$$j(w/5)=j(E_{5}(\tau(b))) = j(E_5(\xi^5)) = \frac{(1-12\xi^5+14\xi^{10}+12\xi^{15}+\xi^{20})^3}{\xi^{25}(1-11\xi^5-\xi^{10})},$$
while $\xi^5=\tau(\eta^5)$ and (2.2) imply that
$$j(w/5)=\frac{(1+228\eta^5+494\eta^{10}-228\eta^{15}+\eta^{20})^3}{\eta^5(1-11\eta^5-\eta^{10})^5}.\eqno{(3.10)}$$
In the same way we have
\begin{align*}
j(w) &= \frac{(1-12\eta^5+14\eta^{10}+12\eta^{15}+\eta^{20})^3}{\eta^{25}(1-11\eta^5-\eta^{10})}\\
&=\frac{(1+228\xi^5+494\xi^{10}-228\xi^{15}+\xi^{20})^3}{\xi^5(1-11\xi^5-\xi^{10})^5}.
\end{align*}
It follows that the minimal polynomials of $\xi$ and $\eta$ divide the polynomial $F_d(x^5)$, where $F_d(x)$ is given by (1.2), as well as the polynomial $G_d(x^5)$, where
$$G_d(x^5)=x^{5h(-d)}(1-11x^5-x^{10})^{5h(-d)} H_{-d}(j_{5,5}(x^5)). \eqno{(3.11)}$$

\section{Fields generated by values of $r(\tau)$.}

If $R_d(X)$ is the minimal polynomial of $z=b-1/b$ over $\mathbb{Q}$, as in Proposition 3.2, define the polynomial $Q_d(X)$ by
$$Q_d(X) = X^{2h(-d)} R_d\left(X-\frac{1}{X} \right).\eqno{(4.1)}$$
The case $d=4$ is unusual, in that
$$F_4(x)=(x^2+1)^2(x^4+18x^3+74x^2-18x+1)^2$$
is divisible by a square factor, so that $Q_4(x)=x^2+1$.  In all other cases we have the following result.  We will need the well-known fact that
$$-z-11=x_1(w)^3 \cong \wp_5'^3.\eqno{(4.2)}$$
(See \cite{deu2}, p.32.)  \bigskip

\noindent {\bf Proposition 4.1.} {\it If $d > 4$, the polynomial $Q_d(x)$ defined by (4.1) is an irreducible factor of $F_d(x)$ of degree $4h(-d)$.  Both $b$ and $\tau(b)$ are roots of $Q_d(x)$.  Furthermore, $Q_d(x^5)$ is divisible by an irreducible factor $p_d(x)$ of degree $4h(-d)$ having $\eta$ as a root.}  \smallskip

\noindent {\it Proof.} Certainly, $b$ is a root of $Q_d(x)$.  If $Q_d(x)$ were reducible, it would have to factor into a product of two polynomials of degree $2h(-d)$ over $\mathbb{Q}$.  Neither of these polynomials would be invariant under $z \rightarrow U(z)=\frac{-1}{z}$, since this would imply that $R_d(x)$ factors.  Hence, $b$ would have to lie in $\Omega_f$, and
$$Q_d(x)=f(x)\cdot x^{2h(-d)}f(-1/x)$$
for some irreducible $f(x)$ having $b$ as a root.  Next, note that
$$\tau(b)-\frac{1}{\tau(b)}=\bar \varepsilon^5 \frac{b-\varepsilon^5}{b-\bar \varepsilon^5}+ \varepsilon^5 \frac{b- \bar \varepsilon^5}{b- \varepsilon^5}=\frac{-11b^2+4b+11}{b^2+11b-1}=\frac{-11z+4}{z+11}.$$
Putting $z_1=\tau(b)-\frac{1}{\tau(b)}$, the last equation gives
$$-z_1-11=\frac{125}{-z-11}=\frac{125}{x_1(w)^3}=x_1(-5/w)^3,$$
by the well-known transformation formula $\eta(-1/\tau)=\sqrt{\frac{\tau}{i}}\eta(\tau)$ for the Dedekind $\eta$-function.  Furthermore,
$$\frac{-5}{w}=\frac{-5w'}{N(w)}=\frac{-w'}{a}=\frac{-v+\sqrt{-d}}{2a}$$
is an ideal basis quotient for the ideal $\mathfrak{a}'=(a,-w')$, where $\wp_5 \mathfrak{a}=(w)$ and therefore $\wp_5' \mathfrak{a'}=(-w')$.  It follows that
$$x_1(-5/w)^3=\left (\frac{\eta\left(\frac{-w'}{5a}\right)}{\eta\left(\frac{-w'}{a}\right)} \right)^6=\overline{x_1(w/a)^3}.$$
From \cite{deu2}, p.32, we have with $z_2 = \bar z_1$ that
$$-z_2-11=x_1(w/a)^3 \cong \wp_5'^3 \cong -z-11$$
and $J(z_2)=j(w/a)$, in the notation of Lemma 2.2.  That lemma implies that $z_2=z^{\sigma^{-1}}$ is a conjugate of $z$ over $K$.  Hence $z_1$ is a conjugate of $z$ over $\mathbb{Q}$, and therefore also a root of $R_d(z)=0$.  This shows that $\tau(b)$ is also a root of $Q_d(x)=0$.  But then either $\tau(b)$ or $\frac{-1}{\tau(b)}$ is a conjugate of $b$ over $\mathbb{Q}$.  From the formula (2.1) for $\tau(b)$, which is linear fractional in $\varepsilon^5$ with determinant $b^2+1 \neq 0$ (for $d > 4$), this would imply that $\sqrt{5} \in \Omega_f$, which is not the case, since $p=5$ is not ramified in $\Omega_f$.  Therefore $Q_d(x)$ is irreducible over $\mathbb{Q}$.  \smallskip

The last assertion of this proposition follows from the equation $\eta^5=b$ and the above arguments.  We have chosen $\eta$ so that $\eta \in \mathbb{Q}(b)$, so the minimal polynomial of $\eta$, namely $p_d(x)$, has degree $4h(-d)$.  $\square$ \bigskip

As a corollary of this argument we have: \bigskip

\noindent {\bf Corollary 4.2.}  {\it The roots of $R_d(z)=0$ are invariant under the map} $z \rightarrow \frac{-11z+4}{z+11}$:
$$(z+11)^{2h(-d)} R_d\left(\frac{-11z+4}{z+11}\right)=5^{3h(-d)}R_d(z).$$ \medskip

Note that the substitution $z \rightarrow V(z)=\frac{-11z+4}{z+11}$ has the effect of  interchanging $j(w)$ and $j(w/5)$, as functions of $z=b-\frac{1}{b}$.  \bigskip

\noindent {\bf Proposition 4.3.} {\it If $d > 4$, the minimal polynomial $p_d(x)$ of $\eta=b^{1/5}$ over $\mathbb{Q}$ is irreducible and normal over $L=\mathbb{Q}(\zeta_5)$.  Furthermore,
$$F=(\Sigma_5 \Omega_f \ or \ \Sigma_5\Omega_{5f})= \mathbb{Q}(b, \zeta_5)=\mathbb{Q}(\eta,\zeta_5)$$
is the disjoint compositum of $\mathbb{Q}(b)=\mathbb{Q}(\eta)$ and $\mathbb{Q}(\zeta_5)$ over $\mathbb{Q}$.  The same facts hold with $b$ replaced by $\tau(b)$ and $\eta$ replaced by $\xi$.}  \smallskip

\noindent {\it Proof.}  We know that a root of $p_d(x)$ generates a quadratic extension of $\Omega_f$ over $\mathbb{Q}$.  Hence, the field $L(\eta)$ contains $L \Omega_f$.  On the other hand, the roots $u$ of (3.1) are contained in $L(\eta)$, since $\xi=(\varepsilon u)^{-1}$ lies in $\mathbb{Q}(\tau(b)) \subseteq \mathbb{Q}(b, \sqrt{5}) \subseteq L(\eta)$, by Theorem 3.3.  Since the $X$-coordinates of points in $E_5[5]$ generate $F$ over $\Omega_f$, and these $X$-coordinates are rational functions in $u$ with coefficients in $L$, by the formulas in  \cite{mor6}, it follows that $F = L(\eta)=\mathbb{Q}(b, \zeta_5)$, and therefore $[L(\eta):L]=\frac{16h(-d)}{4}=4h(-d)$.  This shows that $p_d(x)$ is irreducible over $L=\mathbb{Q}(\zeta_5)$ and implies that $\mathbb{Q}(b) \cap \mathbb{Q}(\zeta_5) = \mathbb{Q}$.
$\square$ \bigskip

This proposition also shows that the polynomial $Q_d(x)$ is not normal over $\mathbb{Q}$, since it has both $b$ and $\tau(b)$ as roots, and $\sqrt{5} \notin \mathbb{Q}(b)$.  Hence, $p_d(x)$ is also not normal over $\mathbb{Q}$.  But $\mathbb{Q}(b) \subset F$ is abelian over $K$ and $\mathbb{Q}(b)$ and $\Omega_f(\zeta_5)$ are linearly disjoint over $\Omega_f$.  \bigskip

\noindent {\bf Corollary 4.4.} {\it If $Q_d(x^5) = p_d(x) q_d(x)$, then $q_d(x)$ is irreducible over $\mathbb{Q}$, of degree $16h(-d)$, and $p_d(\xi)=0$.  Moreover, $x^{4h(-d)}p_d(-1/x)=p_d(x)$ and $x^{16h(-d)}q_d(-1/x)=q_d(x)$.}  \smallskip

\noindent {\it Proof.} To show that the polynomial $q_d(x)$ in $Q_d(x^5) = p_d(x) q_d(x)$ is irreducible, note that $b \in \mathbb{Q}(\zeta \eta)$ implies $\eta$ and therefore also $\zeta$ lies in this field.  Thus, $\mathbb{Q}(\zeta \eta)=\mathbb{Q}(\zeta, \eta)=F$ has degree $8$ over $\Omega_f$ and degree $16h(-d)$ over $\mathbb{Q}$.  This implies that $\zeta \eta$, which is a root of $Q_d(x^5)$, must be a root of $q_d(x)$, hence $q_d(x)$ is irreducible.  Since the set of roots of $Q_d(x^5)$ is stable under the mapping $x \rightarrow -1/x$ and $p_d(x)$ and $q_d(x)$ have different degrees, the respective sets of roots of the latter polynomials must also be stable under this map.  The fact that $x^{4h(-d)}p_d(-1/x)=p_d(x)$ now follows from the norm formula
$$N_{\mathbb{Q}(\eta)/\mathbb{Q}}(\eta) = N_{\Omega_f/\mathbb{Q}}(N_{\mathbb{Q}(\eta)/\Omega_f}(\eta))=1,$$
since $\eta$ is a unit and $\Omega_f$ is complex.  Finally, $\xi$ must also be a root of $p_d(x)$, since $\xi$ and $\tau(b)$ have degree $4h(-d)$ over $\mathbb{Q}$, by Proposition 4.1.  $\square$  \bigskip 

This corollary allows us to prove the following. \bigskip

\noindent {\bf Theorem 4.5.}  {\it The quantities $\eta$ and $\xi$ satisfy
$$\eta=\delta r^{\delta}\left(\frac{w}{5}\right), \ \ \xi = \delta \zeta^{\delta j} r^{\delta}\left(\frac{-1}{w}\right), \ \ \delta=\pm 1, \ \zeta^j \neq 1 \eqno{(4.3)}$$
and are roots of $p_d(x)$.  Thus, the roots of $p_d(x)$ are conjugates over $\mathbb{Q}$ of the values $r(w/5)$ and $\zeta^j r(-1/w)$ of the Rogers-Ramanujan function $r(\tau)$}. \medskip

\noindent {\it Proof.} First note that the map $\sigma: b \rightarrow -1/b$ is an automorphism of $\mathbb{Q}(b)$ which fixes $\Omega_f=\mathbb{Q}(z)$.  Since $\eta$ is the only fifth root of $b$ contained in $\mathbb{Q}(b)$, this automorphism takes $\eta$ to $\eta^\sigma=-1/\eta$ and therefore $\eta-1/\eta \in \Omega_f$.  Furthermore, $\eta'=\zeta \eta$ is a root of the polynomial $q_d(x)$ in Corollary 4.4, and $\eta' \rightarrow -1/\eta'$ is likewise an automorphism of order $2$ of the field $F$.  But then $\eta'-1/\eta'$ has degree $8h(-d)$ over $\mathbb{Q}$, since $\eta'$ is a primitive element for $F$ over $\mathbb{Q}$, so that $\eta'-1/\eta' \notin \Omega_f$.  On the other hand, the function $r(\tau)$ satisfies the identity
$$r^{-1}(\tau)-1-r(\tau) = \frac{\eta(\tau/5)}{\eta(5\tau)},$$
by  \cite{du}, p. 149.  Putting $\tau=w/5$ therefore gives that
$$r(w/5)-r^{-1}(w/5) = -1-\frac{\eta(w/25)}{\eta(w)}=-1-y(w) \in \Omega_f.$$
Now the first formula in (3.9) implies that $i=0$, i.e., that the first formula in (4.3) holds.  On the other hand, putting $\tau=-1/w$ gives
$$r(-1/w)-r^{-1}(-1/w) = \frac{\bar \varepsilon r(w)+1}{r(w)-\bar \varepsilon}-\frac{r(w)-\bar \varepsilon}{\bar \varepsilon r(w)+1}=-\frac{r^2(w)-4r(w)-1}{r^2(w)+r(w)-1},\eqno{(4.4)}$$
and the last expression is linear fractional (with determinant $-5$) in the expression
$$r(w)-r^{-1}(w)=-1-\frac{\eta(w/5)}{\eta(5w)}=-1-y(5w).\eqno{(4.5)}$$
In this case, $y(5w) \in \Omega_{5f}$ ( \cite{sch}, p. 159), but $y(5w) \not \in \Omega_f$, since
$$y(5w)^{24}=\left(\frac{\eta(w/5)}{\eta(w)}\right)^{24} \left(\frac{\eta(w)}{\eta(5w)}\right)^{24}=x_1(w)^{12} \frac{\Delta(w,1)}{\Delta(5w,1)}=x_1(w)^{12} \frac{5^{12}}{\varphi_{P}(w)},$$
where $P$ is the $2 \times 2$ diagonal matrix with entries $5$ and $1$ (in the notation of Hasse  \cite{h} and Deuring  \cite{deu2}).  By  \cite{deu2}, p.43, $\varphi_P(w)$ is a unit, so this gives that $y(5w)^{24} \cong \wp_5'^{12} 5^{12} = \wp_5'^{24} \wp_5^{12}$, i.e. $y(5w)^2 \cong \wp_5'^2 \wp_5$.  This equation implies that $\wp_5$ is the square of an ideal in $\Omega_f(y(5w))$, which shows that $y(5w) \not \in \Omega_f$.  Since $\xi-\xi^{-1} \in \Omega_f$, this proves that $\zeta^j \neq 1$ in (3.9), i.e. that (4.3) holds.  $\square$. \bigskip

\noindent {\bf Theorem 4.6.} {\it If $d \neq 4f^2$ and $z=b-\frac{1}{b}$ is given by (2.5), then $\mathbb{Q}(b)=\Sigma_{\wp_5'} \Omega_f$ is the compositum of $\Omega_f$ with the ray class field of conductor $\wp_5'$ over $K$; and $\mathbb{Q}(\tau(b)) = \Sigma_{\wp_5} \Omega_f$.  Furthermore, the normal closure of $\mathbb{Q}(b)$ over $\mathbb{Q}$ is $\mathbb{Q}(b, \sqrt{5})= \Sigma_{\wp_5} \Sigma_{\wp_5'} \Omega_f$.}  \medskip

\noindent {\it Proof.} First note that $\displaystyle [\Sigma_{\wp_5'}: \Sigma]=\frac{\phi(\wp_5')}{2}=2$, so that $[\Sigma_{\wp_5'}\Omega_f: \Omega_f]=2$.  Moreover, the quadratic extensions $\Sigma_{\wp_5'}\Omega_f$ and $\Sigma_{\wp_5}\Omega_f$ are contained in $F=\Sigma_5 \Omega_f$, because $\Sigma_{\wp_5'}, \Sigma_{\wp_5} \subset \Sigma_5$.  On the other hand, $\textrm{Gal}(F/\Omega_f) \cong \mathbb{Z}_2 \times \mathbb{Z}_4$, so that $F$ has three quadratic subfields over $\Omega_f$.  These subfields are $F_1=\Omega_f(b), F_2=\Omega_f(\tau(b)), F_3=\Omega_f(\sqrt{5})$.  The field $F_3$ is normal over $\mathbb{Q}$, while $F_1$ and $F_2$ must coincide with the fields $\Sigma_{\wp_5'}\Omega_f$ and $\Sigma_{\wp_5} \Omega_f$.  The quantity $b$ satisfies the equation $b^2-bz-1=0$, whose discriminant $z^2+4 = (z+1)(z-1)+5$ is divisible by $\wp_5'$ (by (4.2)).  Now note the congruence (from (1.5))
$$j(w) \equiv -\frac{(z^2+2z+1)^3}{z+1} \equiv -(z+1)^5 \ (\textrm{mod} \ \wp_5).$$
This implies that $j(w)$ is conjugate to $-(z+1)$ (mod $\mathfrak{p}$) for every prime divisor $\mathfrak{p}$ of $\wp_5$ in $\Omega_f$.  Further, the discriminant of $H_{-d}(x)$ is not divisible by $p=5$, since the Legendre symbol $\left(\frac{-d}{5} \right)=+1$ (see \cite{deu}).  Hence, the minimal polynomial $m_d(x)$ of $z$ over $K$ satisfies
$$m_d(z) \equiv (-1)^{h(-d)} H_{-d}(-z-1) \ (\textrm{mod} \ \wp_5),$$
and factors into irreducibles of degree $f_1=\textrm{ord}(\wp_5)$, where $f_1$ is the order of $\wp_5$ in the ring class group (mod $f$) of $K$.  If $f_1 \ge 2$, then certainly $z=1$ is not a root of $m_d(z)$ (mod $\wp_5$), so no prime divisor of $\wp_5$ divides $z-1$.  If $f_1=1$, then by the calculations of Proposition 3.2, $d$ is $11$ or $19$ (since $d \neq 16$ by assumption); and it can be checked that
$$R_{11}(z) \equiv (z+1)(z+3), \ \ R_{19}(z) \equiv z(z+1) \ \ (\textrm{mod} \ 5).$$
It follows that no prime divisor of $\wp_5$ divides $z-1$, for any $d$.  Hence, only the prime divisors of $\wp_5'$ in $\Omega_f$ can be ramified in $\Omega_f(b)/\Omega_f$. It follows that $\wp_5'$ must divide the conductor of $F_1$, which proves the first assertion.  Then the field $\Sigma_{\wp_5} \Sigma_{\wp_5'} \Omega_f=F_1F_2$ is obviously the smallest normal extension of $\mathbb{Q}$ containing $ \mathbb{Q}(b)$.  $\square$  \medskip

\noindent {\bf Corollary 4.7.} {\it If $d \neq 4f^2$, $w$ is defined by (2.4), and $\zeta^j$ is as in (4.3), then}
$$\mathbb{Q}(r(w/5)) = \mathbb{Q}(b)=\Sigma_{\wp_5'} \Omega_f, \ \ \mathbb{Q}(\zeta^j r(-1/w)) = \mathbb{Q}(\tau(b))=\Sigma_{\wp_5} \Omega_f,$$
{\it and $\mathbb{Q}(r(-1/w))=\mathbb{Q}(r(w)) = F =\Sigma_5 \Omega_f$.  The field $F_1=\mathbb{Q}(\eta)=\mathbb{Q}(r(w/5))$ is the inertia field for $\wp_5$ in the abelian extension $F/K$.} \medskip

\noindent {\it Proof.} The first assertion follows directly from Theorems 4.5 and 4.6, since $\mathbb{Q}(r(w/5))=\mathbb{Q}(\eta)=\mathbb{Q}(b)$.  The fact that $\mathbb{Q}(r(-1/w))= F$ follows from $r^\delta(-1/w)=\delta \zeta^{-\delta j}\xi$ and the proof of Corollary 4.4, which shows that $\zeta^{-\delta j} \xi$ is a root of the irreducible polynomial $q_d(x)$.  Now, $r(w)$ generates a field over $\mathbb{Q}$ whose degree is at least $4h(-d)$, by (4.5), since $r(w)$ lies in a quadratic extension of $\Omega_{5f} \subset F$.  If $[\mathbb{Q}(r(w)):\mathbb{Q}]=4h(-d)$, then (4.4) shows that $r(-1/w)$ would have degree at most $8h(-d)$ over $\mathbb{Q}$, which is not the case, as we have just shown.  Hence, $r(w)$ must have degree at least $4$ over $\Omega_f$.  If this degree equals $4$, then $\mathbb{Q}(r(w))/\Omega_f \subseteq F/\Omega_f$ is a quartic extension which contains $\sqrt{5}$.  (This is easiest to see by using the correspondence between abelian extensions of $\Omega_f$ and characters of $\textrm{Gal}(F/\Omega_f)$, as in  \cite{h3}, p. 5.)  Therefore $r(-1/w) \in \mathbb{Q}(r(w))$ by the first equality in (4.4).  This contradiction proves that $r(w)$ has degree $16h(-d)$ over $\mathbb{Q}$ and $\mathbb{Q}(r(w)) = F$.  The last assertion follows from the fact that the ramification index of the prime divisors of $\wp_5$ in $F/K$ is $e=4=[F:F_1]$, so that $F_1$ is the maximal subextension of $F$ which is unramified at $\wp_5$. 
$\square$ \bigskip

In the case $K=\mathbb{Q}(i)$, we have $\Sigma_{\wp_5}=\Sigma_{\wp_5'}=K$, so the conclusion of Theorem 4.6 cannot hold.  However, the fact that $\wp_5'$ ramifies and $\wp_5$ does not ramify in the quadratic extension $\Omega_f(b)/\Omega_f$ follows in exactly the same way, since $R_{16}(z) \equiv (z+1)(z+2)$ (mod $5$).  This gives the following result.  \bigskip

\noindent {\bf Theorem 4.8.} {\it If $K=\mathbb{Q}(i)$, $d=4f^2>4$ and $2 \mid f$, then with the value of $j$ in (4.3),}
$$\mathbb{Q}(r(w/5)) = \mathbb{Q}(b)=\Sigma_{2\wp_5'}\Omega_f \ \ and \ \ \mathbb{Q}(\zeta^j r(-1/w)) = \mathbb{Q}(\tau(b)) = \Sigma_{2\wp_5} \Omega_f.$$
{\it In general, if $d=4f^2 >4$, then $\mathbb{Q}(r(-1/w))=\mathbb{Q}(r(w)) = F =\Sigma_5 \Omega_{5f}$; and $F_1=\mathbb{Q}(\eta)$ is the inertia field for $\wp_5$ in the abelian extension $F/K$.} \medskip

\noindent {\it Proof.}  In this case we have $f=2f'$ and $\Omega_{5f} = \Omega_{10} \Omega_f$, by Hasse's Zusatz in \cite{h2}, p. 326.  Therefore $F=\Sigma_5 \Omega_{10} \Omega_f$.  On the other hand, $\textsf{S}_5 \cap \textsf{P}_{10} \subset \textsf{S}_{2 \wp_5'}$ in $K=\mathbb{Q}(i)$, when these ideal groups are declared modulo $10$, so we have that $\Sigma_{2 \wp_5'} \subset \Sigma_5 \Omega_{10}$ and $\Sigma_{2 \wp_5'} \Omega_f \subset F$.  Since $[\Sigma_{2 \wp_5'}:K]=2$ and $\wp_5'$ ramifies in $\Sigma_{2 \wp_5'}$, it is clear that $[\Sigma_{2 \wp_5'} \Omega_f :\Omega_f]=2$.  Now the proof of Theorem 4.6 shows that $\mathbb{Q}(b) = \Sigma_{2 \wp_5'} \Omega_f$ and $\mathbb{Q}(\tau(b)) = \Sigma_{2 \wp_5} \Omega_f$ and the rest is a consequence of Theorem 4.5 and the same arguments as in the previous corollary.  $\square$  \bigskip

\noindent {\bf Remark.} When $K=\mathbb{Q}(i)$ and $f$ is odd, the conductor $\mathfrak{f}(F_1/K)$ of $F_1/K$ divides $\wp_5'(f)$, and is divisible by the conductor $\mathfrak{f}(\Omega_f/K)$.  Since $f$ is odd, $\mathfrak{f}(\Omega_f/K)=(f)$, so that $\mathfrak{f}(F_1/K)=\wp_5'(f)$.  (See \cite{co}, Ex. 9.20, pp. 195-196.)  In the general case $d>4$ it is not hard to see that the equality $\mathfrak{f}(F_1/K)=\wp_5'(f)$ still holds, unless $-d =d_K f^2 \neq -4f^2$, $d_K \equiv 1$ (mod $8$), and $f=2f'$ with odd $f'$; in which case $\mathfrak{f}(F_1/K)=\wp_5'(f')$.  As an example of the latter phenomenon, see the polynomial $p_{124}(x)$ in Table 2 below, for which $f=2$, but whose discriminant is not divisible by $2$.  \medskip

Table 1 gives the minimal polynomials $p_d(x)$ of the values $r(w/5)$ for all $d < 150$.  For most values of $d$, $p_d(x)$ was computed from $H_{-d}(x)$ using the fact that $p_d(x) \mid F_d(x^5)$ with $F_d(x)$ in (1.2).  For $d \neq 4f^2$ for which $H_{-d}(x)$ was not available, $p_d(x)$ was computed by approximating to high accuracy the values of $r(\tau) = r(w/(5a))$ at ideal basis quotients of representatives $\wp_5 \mathfrak{a}=(5a, w)$ of the classes in the {\it ray class group modulo} $\mathfrak{f}=\wp_5'$ of $\textsf{R}_{-d}$, for which $\wp_5^2 \mid (w)$, in line with (2.4).  (See \cite{sch}, p.88.)  This gives $2h(-d)$ values $r(w/(5a))$, which are class invariants for the ideal class group $\textsf{A}/\textsf{H}_{\wp_5'f}$, where $\textsf{A}$ is the group of fractional ideals of $K$ prime to $\wp_5' (f)$ and $\textsf{H}=\textsf{H}_{\wp_5'f}$ is the ideal group of conductor $\wp_5' (f)$ (or $\wp_5'(f')$) corresponding to the class field $\mathbb{Q}(r(w/5))/K$.  Then
$$p_d(x) = \prod_{\mathfrak{a} \ mod \ \textsf{H}}{(x-r\left(\frac{w}{5a}\right))(x-\bar r\left(\frac{w}{5a}\right))}.$$
A similar computation was carried out for $d=4f^2$.  In Section 5 below we will give an algebraic method for verifying these calculations.  The discriminants of these polynomials seem to satisfy the following. \bigskip

\noindent {\bf Conjecture.} \begin{enumerate}[(a)]
\item{{\it If $q > 5$ is a prime which divides $d_K$ but does not divide $f$, then $q^{2h(-d)}$ exactly divides $\textrm{disc}(p_d(x))$.}}
\item{{\it If $h=h(-d)$, $5^{h(2h-1)}$ exactly divides $\textrm{disc}(p_d(x))$.}}
\item{{\it $\textrm{disc}(p_d(x))$ is only divisible by primes $q \le d$.}}
\item{{\it If $q \neq 5$ is a prime dividing $\textrm{disc}(p_d(x))$, then the Kronecker symbol $\left(\frac{-d}{q}\right) \neq 1$.}}
\end{enumerate}

\begin{table}
  \centering 
  \caption{ The minimal polynomial $p_d(x)$ of $r(w/5)$, $w = \frac{v+\sqrt{-d}}{2}, \ 5^2 \mid N(w)$, \ $11 \le d \le 99$.}\label{ }

\noindent \begin{tabular}{|c|l|c|}
\hline
  &  &  \\
$d$	&   $\ \ \ p_d(x)$  &   $\textrm{disc}(p_d(x))$ \\
\hline
  &   &   \\
 11  &  $x^4-x^3+x^2+x+1$  &  $5 \cdot 11^2$  \\
 16  &  $x^4-2x^3+2x+1$      &   $2^6 5$    \\
 19 &   $x^4+x^3+3x^2-x+1$  &   $5 \cdot19^2$ \\
  24 &  $x^8-2x^7+x^6-4x^5+3x^4+4x^3+x^2+2x+1$  &  $2^{12} 3^4 5^6$  \\
  31 & $x^{12}-x^{11}+5x^{10}-4x^9+8x^8-2x^7+19x^6+2x^5+8x^4$ & $3^8 5^{15} 31^6$\\
  & $+4x^3+5x^2+x+1$ & \\
  36  &  $x^8 +x^6-6x^5+9x^4+6x^3+x^2+1$  &  $2^8 3^6 5^6 11^4$ \\
  39 & $x^{16}-3x^{15}+7x^{14}-9x^{13}+21x^{12}-15x^{11}+17x^{10}+3x^9$ & $3^8 5^{28} 7^8 13^8$\\  
  & $+11x^8-3x^7+17x^6+15x^5+21x^4+9x^3+7x^2+3x+1$ & \\
  44 & $x^{12}-x^{11}+6x^{10}+15x^8+9x^6+15x^4+6x^2+x+1$ & $2^8 5^{15} 11^6 19^4$\\
  51  &   $x^8+x^7+x^6-7x^5+12x^4+7x^3+x^2-x+1$  &  $2^{12} 3^4 5^6 17^4$  \\
  56 & $x^{16}+8x^{14}-4x^{13}+15x^{12}-12x^{11}+50x^{10}+4x^9+91x^8$ & $2^{40}5^{28}7^831^4$\\
   & $-4x^7+50x^6+12x^5+15x^4+4x^3+8x^2+1$ & \\
  59 & $x^{12}-4x^{11}+5x^{10}-2x^9+14x^8-2x^7-24x^6+2x^5+14x^4$ & $2^{20} 5^{15} 59^6$ \\
   & $+2x^3+5x^2+4x+1$  & \\
  64 & $x^8+4x^7+10x^6+8x^5+12x^4-8x^3+10x^2-4x+1$ & $2^{18} 3^8 5^6$ \\
   71 & $x^{28}-6x^{27}+17x^{26}-45x^{25}+104x^{24}-164x^{23}+277x^{22}$ & $5^{91} 7^{16} 23^8 71^{14}$\\
     & $-357x^{21}+388x^{20}-319x^{19}+316x^{18}+135x^{17}-144x^{16}+83x^{15}$ & \\
     & $-551x^{14}-83x^{13}-144x^{12}-135x^{11}+316x^{10}+319x^9+388x^8$ & \\
     & $+357x^7+277x^6+164x^5+104x^4+45x^3+17x^2+6x+1$ & \\
    76 & $x^{12}-5x^{11}+12x^{10}-2x^9-21x^8+12x^7+35x^6-12x^5$ & $2^8 3^{12} 5^{15} 19^6$\\
   & $-21x^4+2x^3+12x^2+5x+1$ & \\
   79 & $x^{20}+9x^{18}-12x^{17}+18x^{16}-9x^{15}+117x^{14}-33x^{13}+99x^{12}$ & $3^{28} 5^{45} 29^8 79^{10}$\\ 
   & $-207x^{11}+353x^{10}+207x^9+99x^8+33x^7+117x^6+9x^5$ & \\
   & $+18x^4+12x^3+9x^2+1$ & \\
 84  &  $x^{16} + 2x^{15} -4x^{14} -12x^{13} + 25x^{12}-18x^{11}+68x^{10}-112x^9$ &  $2^{32}3^{20} 5^{28}7^859^4$\\
  & $+13x^8+112x^7+68x^6+18x^5+25x^4+12x^3-4x^2-2x+1$  &  \\
 91  &  $x^8 +4x^7 -x^6-14x^5+23x^4+14x^3-x^2-4x+1$  &  $2^8 3^4 5^6 7^4 13^4$  \\
 96  &  $x^{16} + 4x^{15} + 29x^{12}-24x^{11}+86x^{10}-32x^9+105x^8 +32x^7$ &  $2^{32}3^{24} 5^{28}71^4$\\
  & $+86x^6+24x^5+29x^4-4x+1$  &  \\
 99  & $x^8+7x^7+15x^6+15x^5+16x^4-15x^3+15x^2-7x+1$  & $2^{12}3^4 5^6 11^4$  \\   
  &  &  \\
\hline
\end{tabular}
\end{table}

\begin{table}
  \centering 
  \caption{ The minimal polynomial $p_d(x)$ of $r(w/5)$, $w = \frac{v+\sqrt{-d}}{2}, \ 5^2 \mid N(w)$,  \ $104 \le d \le 144$.}\label{ }
\noindent \begin{tabular}{|c|l|c|}
\hline
  &  &  \\
$d$	&   $\ \ \ p_d(x)$  &   $\textrm{disc}(p_d(x))$ \\
\hline
  &   &   \\
  104 & $x^{24}-4x^{23}+20x^{22}-40x^{21}+53x^{20}-28x^{19}+94x^{18}-92x^{17}$ & $2^{84} 5^{66} 13^{12} 29^8 79^4$\\
   & $+42x^6-76x^{15}+782x^{14}-328x^{13}-272x^{12}+328x^{11}+782x^{10}$ & \\
   & $76x^9+42x^8+92x^7+94x^6+28x^5+53x^4+40x^3+20x^2$ & \\
   & $+4x+1$ & \\
   111 & $x^{32}-4x^{31}+21x^{30}-31x^{29}+144x^{28}-180x^{27}+563x^{26} $ & $3^{52} 5^{120} 11^{12} 37^{16}$\\
   & $-435x^{25}+1398x^{24}-653x^{23}+2108x^{22}+380x^{21}+4093x^{20}$ & $\times 43^8 61^8$\\
   & $+1273x^{19}+4560x^{18}-990x^{17}+7975x^{16}+990x^{15}+4560x^{14}$ & \\
   & $-1273x^{13}+4093x^{12}-380x^{11}+2108x^{10}+653x^9+1398x^8$ & \\
   & $+435x^7+563x^6+180x^5+144x^4+31x^3+21x^2+4x+1$ & \\
  116 & $x^{24}-6x^{23}+12x^{22}-24x^{21}+99x^{20}-58x^{19}+136x^{18}-256x^{17}$ & $2^{80}5^{66}7^8 29^{12} 41^8$ \\
 & $+144x^{16}+410x^{15}+436x^{14}+274x^{13}-1192x^{12}-274x^{11}$ & \\
 & $+436x^{10}-410x^9+144x^8+256x^7+136x^6+58x^5+99x^4$ & \\
 & $+24x^3+12x^2+6x+1$ & \\
 119 & $x^{40}-x^{39}+12x^{38}-51x^{37}+146x^{36}-248x^{35}+569x^{34}$ & $5^{190} 7^{20} 11^{24} 17^{20}$ \\
  & $-951x^{33}+2005x^{32}-3810x^{31}+8702x^{30}-14440x^{29}$ & $\times 19^{12} 23^{16} 47^8$ \\
  & $+26580x^{28}-35295x^{27}+47491x^{26}-45351x^{25}+53426x^{24}$ & \\
  & $-29809x^{23}+41387x^{22}-6812x^{21}+31769x^{20}+6812x^{19}$ & \\
  & $+41387x^{18}+29809x^{17}+53426x^{16}+45351x^{15}+47491x^{14}$ & \\
  & $+35295x^{13}+26580x^{12}+14440x^{11}+8702x^{10}+3810x^9$ & \\
  & $+2005x^8+951x^7+569x^6+248x^5+146x^4+51x^3$ & \\
  & $+12x^2+x+1$ & \\
  124 & $x^{12}-7x^{11}+9x^{10}+8x^9+24x^8+6x^7-67x^6-6x^5+24x^4$ & $3^{12}5^{15} 11^4 31^6$\\
   & $-8x^3+9x^2+7x+1$ & \\
  131 & $x^{20}+20x^{18}+8x^{17}+48x^{16}+4x^{15}+72x^{14}+88x^{13}+348x^{12}$ & $2^{76} 5^{45} 31^4 131^{10}$\\
   & $168x^{11}+446x^{10}-168x^9+348x^8-88x^7+72x^6-4x^5+48x^4$ & \\
   & $-8x^3+20x^2+1$ & \\
  136 & $x^{16}+6x^{15}+25x^{14}+24x^{13}-3x^{12}+119x^{10}+174x^9+404x^8$ & $2^{56} 3^{16} 5^{28} 11^8 17^8$\\
    & $-174x^7+119x^6-3x^4-24x^3+25x^2-6x+1$ & \\
 139 & $x^{12}-5x^{11}+12x^{10}+16x^9+33x^8+12x^7-55x^6-12x^5$ & $2^{24} 3^{12} 5^{15} 139^6$ \\
  & $+33x^4-16x^3+12x^2+5x+1$ & \\
  144 & $x^{16}-2x^{15}+18x^{14}+24x^{13}+83x^{12}+78x^{11}+74x^{10}+40x^9$ & $2^{24} 3^{12} 5^{28} 7^8$\\
    & $+9x^8-40x^7+74x^6-78x^5+83x^4-24x^3+18x^2+2x+1$ & $\times 11^4 19^8$ \\
   &  &  \\
\hline
\end{tabular}
\end{table}

\section{Periodic points of an algebraic function.}

\subsection{Preliminary facts on the group $G_{60}$.}

In this section we shall make use of the fact that the rational function
$$f_5(z)=\frac{(1+228z^5+494z^{10}-228z^{15}+z^{20})^3}{z^5(1-11z^5-z^{10})^5}$$
is invariant under a group $G_{60}$ of linear fractional substitutions:
$$G_{60}= \langle S,T \rangle, \ \ S(z)=\zeta z, \ \ T(z)=\frac{-(1+\sqrt{5})z+2}{2z+1+\sqrt{5}},$$
which is isomorphic to the icosahedral group $A_5$.  The coefficients of the maps in $G_{60}$ are in the field $\mathbb{Q}(\zeta_5)$.  The transformations $S$ and $T$ have orders 5 and 2, respectively, while the transformation
$$U(z)=\frac{-1}{z}$$
is given in terms of $S$ and $T$ by $U=T \cdot S^2 \cdot T \cdot S^3 \cdot T \cdot S^2$.  (See \cite{fr1}, II, pp. 42-43.)  Furthermore,
$$H=\{1, T, U, TU\}$$
is a Klein-4 subgroup of $G_{60}$, where $TU(z)=UT(z)=-1/T(z)=T_2(z)$, and
$$T_2(z)=\frac{-(1-\sqrt{5})z+2}{2z+1-\sqrt{5}}.$$
Thus, $U=TT_2=T_2T$.  The normalizer of $H$ in $G_{60}$ is $N=\langle A, H \rangle \cong A_4$, where $A=STS^{-2}$ is the map
$$A(z)=\zeta^3 \frac{(1+\zeta)z+1}{z-1-\zeta^4}$$
of order $3$, and $ATA^{-1}=U, AUA^{-1}=T_2$.  Also, $A^\sigma=A^{-1}U$ is the conjugate map
$$A^\sigma(z)=\zeta \frac{(1+\zeta^2)z+1}{z-1-\zeta^3},$$
obtained by applying the automorphism $\sigma: \zeta \rightarrow \zeta^2$ to the coefficients.  In particular, $\textrm{Gal}(\mathbb{Q}(\zeta)/\mathbb{Q})$ is a subgroup of the automorphism group $\textrm{Aut}(N)$.  \medskip

It is clear from (2.2) and (3.11) that $\textrm{deg}(G_d(x^5))=60h(-d)$.  The group $G_{60}$ acts on the irreducible factors $p(x)$ of $G_d(x^5)$ over $L=\mathbb{Q}(\zeta_5)$, one of which is $p_d(x)$ (Proposition 4.3), by
$$p^\sigma(x) =(cx+d)^{deg(p)}p(\sigma(x)) = (cx+d)^{deg(p)}p\left(\frac{ax+b}{cx+d}\right), \ \ \sigma \in G_{60},$$
ignoring constant factors.  Moreover, $G_{60}$ acts transitively on these irreducible factors over the field $L$ (see the analogous argument in  \cite{lym}, p. 1982), so $G_d(x^5)$ splits into $15$ irreducible factors of degree $4h(-d)$ over $L$.  In particular, these considerations show that every root of $G_d(x^5)$ has the form $\sigma(\alpha)$ for some root $\alpha$ of $p_d(x)$ and some $\sigma \in G_{60}$.  \medskip

The group $G_{60} \cong A_5$ has no elements of order $4$, so the stabilizer of $p_d(x)$ is one of the five conjugate subgroups in $G_{60}$ of the subgroup $H$.  We have that
$$S^{-1}US(z) = \frac{-\zeta^3}{z}, \ \ S^{-1}TS(z)=\frac{-(1+\sqrt{5})z+2\zeta^4}{2\zeta z+(1+\sqrt{5})}.$$
Hence, only one these conjugate subgroups, namely $H$, contains the map $U$, and since $U$ fixes $p_d(x)$ by Corollary 4.4, we have
$$\textrm{Stab}_{G_{60}}(p_d(x)) = H =  \{1, T, U, TU\}.$$
As a consequence, we have that
$$\left(z+\frac{1+\sqrt{5}}{2}\right)^{4h(-d)}p_d(T(z))= \left(\frac{5+\sqrt{5}}{2} \right)^{2h(-d)} p_d(z).$$
It can be checked that the factor on the right side of this equation is correct by putting $z$ equal to
$$z_1=\frac{-1-\sqrt{5}+\sqrt{10+2\sqrt{5}}}{2},$$
which is a fixed point of $T(z)$, and noting that $p_d(z_1) \neq 0$, since $\mathbb{Q}(z_1)$ is a cyclic quartic extension of $\mathbb{Q}$ in which $p=5$ is totally ramified. \medskip

We also note that all of the roots of $p_d(x)$ are values of the Rogers-Ramanujan function $r(\tau)$.  This follows from the identity (see \cite{du}, p. 138):
$$j(\tau) = \frac{(r^{20}-228r^{15}+494r^{10}+228r^{5}+1)^3}{r^5(1-11r^5-r^{10})^5}, \ \ r = r(\tau).$$
Any root $\alpha$ of $p_d(x)$ satisfies $f_5(\alpha) = j(w/a)$ for some $w$ of the form (2.4) and some $a$.  However, equations (2.2), (2.6), and (2.7) imply that $f_5(r(w/a))=j(w/a)$.  It follows that $\alpha$ and $r(w/a)$ are related by an element $M$ of the group $G_{60}$.  Now we use Proposition 2 of \cite{du}, according to which
$$r(\tau+1) = S(r(\tau)), \ \ r\left(\frac{-1}{\tau}\right) = T(r(\tau)) \ \ \tau \in \mathbb{H}.$$
It follows that the action of any mapping $M \in G_{60}$ on a value $r(\tau)$ can be represented by a suitable element $\mu \in \Gamma=SL_2(\mathbb{Z})$, such that $M(r(\tau))=r(\mu(\tau))$; hence,
$$\alpha = M(r(w/a)) = r(\mu(w/a))$$
is a value of the function $r(\tau)$ with $\tau \in K$.  This argument applies to all the roots of $G_d(x^5)$.  (Since $r(\tau)$ is a Hauptmodul for $\Gamma(5)$, the above formulas imply that $G_{60} \cong \bar \Gamma(5)$; see \cite{scho}, p. 76.)
\bigskip

\subsection{Automorphisms of $F_1/K$.}

Now let $\psi$ be an automorphism of the extension $F=\Omega_f(\xi, \zeta_5)$ which fixes $\Omega_f(\xi)=\Omega_f(\tau(b))$ and sends $\zeta$ to $\zeta^2$.  Then $\psi$ takes $\sqrt{5}$ to $-\sqrt{5}$, so that
$$(\eta^5)^{\psi}=b^{\psi}=\tau(\xi^5)^{\psi}=\frac{-\xi^5+\bar \varepsilon^5}{\bar \varepsilon^5 \xi^5+1}=-\frac{\varepsilon^5 \xi^5+1}{-\xi^5+\varepsilon^5}=\frac{-1}{\eta^5}.$$
It follows that $\eta^\psi=\frac{-\zeta^i}{\eta}$, for some $i$.  Thus, $\zeta^i \in \Omega_f(\eta)$ and $i \equiv 0$ (mod $5$),  giving $\eta^\psi=\frac{-1}{\eta}$.   \medskip

Next, let $\phi$ be an automorphism of $F$ which takes $\eta$ to $\xi$ and fixes $\zeta$ (this exists by Proposition 4.3 and Corollary 4.4).  Then
$$\tau(b)^\phi = (\xi^5)^\phi=\tau(\eta^5)^\phi=\tau(\xi^5)=\eta^5=b,$$
so that $\xi^\phi = \eta$ by Theorem 3.3, since $\zeta \notin \mathbb{Q}(b)$.  Hence $\phi$ has order $2$ in $\textrm{Gal}(F/\mathbb{Q})$.  Furthermore, since
$$-z^\phi-11=-\left(b-\frac{1}{b}\right)^\phi-11=-\left(\tau(b)-\frac{1}{\tau(b)}\right)-11=-z_1-11,$$
we see from (4.2) and $-z_1-11 \cong \wp_5^3$ (see the proof of Proposition 4.1) that $\phi$ interchanges the ideals $\wp_5'$ and $\wp_5$.  Thus, $\phi$ does not fix the field $K$.  \medskip

Since $T \in H$, the map $\sigma_1=\left(\eta \rightarrow T(\eta) \right)$ also represents an automorphism of order $2$ of $F/L$.
Setting $\upsilon=\eta-\frac{1}{\eta} \in \Omega_f$, and noting that $\upsilon$ is an algebraic integer, we have
$$T(\eta)-\frac{1}{T(\eta)}=-\frac{\eta^2-4\eta-1}{\eta^2+\eta-1}=-\frac{\upsilon-4}{\upsilon+1}=-1+\frac{5}{\upsilon+1},$$
so that
$$(\upsilon+1)^{\sigma_1}=\frac{5}{\upsilon+1}.\eqno{(5.1)}$$
The identity
$$x^5-\frac{1}{x^5}=\left(x-\frac{1}{x}\right)^5+5\left(x-\frac{1}{x}\right)^3+5\left(x-\frac{1}{x}\right)$$
gives that
$$z=b-\frac{1}{b}=\upsilon^5+5\upsilon^3+5\upsilon,$$
and implies
$$z \equiv \upsilon^5 \ \ (\textrm{mod} \ 5).$$
It follows that
$$z+11 \equiv z+1 \equiv (\upsilon+1)^5 \ \ (\textrm{mod} \ 5),$$
so $\upsilon+1$ is divisible by $\wp_5'$ but not by any prime divisors of $\wp_5$.  Equation (5.1) implies that $(\upsilon + 1) = \left(\frac{\eta^2+\eta-1}{\eta}\right)=\wp_5'$, and that $\sigma_1$ interchanges the ideals $\wp_5$ and $\wp_5'$.  This also shows that
$$\wp_5=\left(\frac{5\eta}{\eta^2+\eta-1}\right)=\left(\frac{\xi^2+\xi-1}{\xi}\right) \ \ \textrm{in} \ \Omega_f.$$

\subsection{Periodic points.}

Thus, the automorphism $\sigma_1 \phi$ fixes the field $K$, and it follows from (3.10) and the fact that $\sigma_1$ fixes the rational function $f_5(\eta)$ that 
$$j(w/5)^{\sigma_1 \phi}=\frac{(1+228\xi^5+494\xi^{10}-228\xi^{15}+\xi^{20})^3}{\xi^5(1-11\xi^5-\xi^{10})^5}=j(w).$$
Since $\sigma_1 \phi$ fixes the quadratic field $K$ and $K(j(w))=\Omega_f$, we deduce that
$$(\sigma_1 \phi) |_{\Omega_f} = \left(\frac{\Omega_f/K}{\wp_5} \right).$$
We would like to extend this automorphism to the abelian extension $F_1=\mathbb{Q}(\eta)=\Omega_f(\eta)$ of $K$, in which $\wp_5$ is still unramified.  This can be done in two ways.  On the one hand, the restriction of
$$ \tau_5=\left(\frac{F_1/K}{\wp_5} \right)= \left(\frac{\mathbb{Q}(b)/K}{\wp_5} \right)$$
to $\Omega_f$ is certainly the same as $(\sigma_1 \phi) |_{\Omega_f}$.  But the automorphism $\rho=\psi |_{F_1}=(\eta \rightarrow \frac{-1}{\eta})$ of $F_1$ fixes $\Omega_f$, so that $\rho \tau_5=\tau_5 \rho \in \textrm{Gal}(F_1/K)$ also restricts to $(\sigma_1 \phi) |_{\Omega_f}$.  Hence we have that
$$ \tau_5= \sigma_1 \phi \ \ \textrm{or} \ \ \tau_5 \rho= \sigma_1 \phi \ \ \textrm{on} \ F_1.$$
This gives
$$\eta^{\tau_5}=\eta^{\sigma_1 \phi}=T(\eta)^\phi=T(\xi), \ \ \textrm{or} \ \ \eta^{\tau_5 \rho}=\eta^{\sigma_1 \phi}=T(\xi).$$
Hence,
$$\xi=T(\eta^{\tau_5}) = \frac{-(1+\sqrt{5})\eta^{\tau_5}+2}{2\eta^{\tau_5}+1+\sqrt{5}} \ \ \textrm{or} \ \ \xi=T_2(\eta^{\tau_5}) = \frac{-(1-\sqrt{5})\eta^{\tau_5}+2}{2\eta^{\tau_5}+1-\sqrt{5}}.$$
In the following theorem we eliminate the second of these possibilities.
 \bigskip

\noindent {\bf Theorem 5.1.} {\it If $\displaystyle \tau_5=\left(\frac{\Omega_f(\eta)/K}{\wp_5} \right)$, the coordinates of the solution $(\xi,\eta)$ of $\mathcal{C}_5$ satisfy
$$\xi= T(\eta^{\tau_5})=\frac{-(1+\sqrt{5})\eta^{\tau_5}+2}{2\eta^{\tau_5}+1+\sqrt{5}}.\eqno{(5.2)}$$}
\medskip

\noindent {\it Proof.}  Assume that $d>4$.  It suffices to show that $T(\xi)=\eta^{\tau_5}$, and to do this we show that $T(\xi) \equiv \eta^5$ (mod $\wp_5$) in $F_1=\mathbb{Q}(\eta)$.  We have
\begin{align*}
T(\xi)-\eta^5  &= T(\xi)-\tau(\xi^5) =\frac{\bar \varepsilon \xi+1}{\xi-\bar \varepsilon}-\frac{-\xi^5+\varepsilon^5}{\varepsilon^5 \xi^5+1}.\\
&= \frac{-\xi+\varepsilon}{\varepsilon \xi+1}+\frac{\xi^5-\varepsilon^5}{\varepsilon^5 \xi^5+1}\\
&= \frac{(5+2\sqrt{5})(\xi^2+1)(\xi-\varepsilon)^2}{(\xi^2+\xi + \frac{3+\sqrt{5}}{2})(\xi^2-\frac{3+\sqrt{5}}{2} \xi+\frac{3+\sqrt{5}}{2})},
\end{align*}
by factoring this rational function in $\xi$ and $\sqrt{5}$ on Maple.  Now multiply this expression by
$$(T(\xi)-\eta^5)^\psi=T_2(\xi)+\frac{1}{\eta^5}.$$
This yields the equation
$$(T(\xi)-\eta^5)\left(T_2(\xi)+\frac{1}{\eta^5}\right)=\frac{5(\xi^2+1)^2(\xi^2+\xi-1)^2}{p_1(\xi)p_2(\xi)}\eqno{(5.3)}$$
in $F_1$, where
$$p_1(\xi) = \xi^4+2\xi^3+4\xi^2+3\xi+1, \ \ p_2(\xi)= \xi^4-3\xi^3+4\xi^2-2\xi+1.$$
Expanding the element $\xi^{-4} p_1(\xi)p_2(\xi)$ of $\Omega_f$ $\pi$-adically in terms of the generating element $\pi=(\xi^2+\xi-1)/\xi$ of $\wp_5$ gives
$$\xi^{-4} p_1(\xi)p_2(\xi)=\pi^4-5\pi^3+15\pi^2-25\pi+25, \ \ \pi = \frac{\xi^2+\xi-1}{\xi},$$
and shows that the squares of prime divisors $\mathfrak{q}$ of $\wp_5$ in $F_1$ exactly divide $p_1(\xi)p_2(\xi)$ (recall that $\wp_5$ is unramified in $F_1$).  This shows that $\frac{(\xi^2+1)^2(\xi^2+\xi-1)^2}{p_1(\xi)p_2(\xi)}$ is a $\mathfrak{q}$-adic integer of $F_1$ for each $\mathfrak{q} \mid \wp_5$, and (5.3) gives that
$$(T(\xi)-\eta^5)\left(T_2(\xi)+\frac{1}{\eta^5}\right)\equiv 0 \ \textrm{mod} \ \wp_5.$$
It follows that $T(\xi) \equiv \eta^5$ or $T_2(\xi) =\frac{-1}{T(\xi)} \equiv \frac{-1}{\eta^5}$ (mod $\mathfrak{q}$) for each $\mathfrak{q}$.  Since $T(\xi)$ and $\eta$ are units, the latter congruence implies that $T(\xi) \equiv \eta^5$ (mod $\mathfrak{q}$), which therefore holds for all $\mathfrak{q}$ dividing $\wp_5$.  Thus we have $T(\xi) \equiv \eta^5$ (mod $\wp_5$).  This implies finally that $T(\xi) = \eta^{\tau_5}$, since $T(\xi)=\eta^{\tau_5 \rho}$ would give $\eta^\rho \equiv \eta$ (mod $\mathfrak {q}$), so $\eta \equiv \pm 2$ (mod  $\mathfrak {q}$) and $z \equiv \pm 1$ (mod $N_{F_1/\Omega_f}(\mathfrak{q})$).  As in the proof of Theorem 4.6, this can only happen when $
f_1=\textrm{ord}(\wp_5)=1$ in the ring class group (mod $f$) of $K$ and $d=11,16,19$.  In these cases $[\mathbb{Q}(\eta):K] =2$, so $\textrm{Gal}(\mathbb{Q}(\eta)/K)=\{1, \rho \}$.  In the first two cases $\tau_5$ has order $2$, so $\tau_5=\rho$, while in the third case $\tau_5=1$.  In all three cases the formula (5.2) can be checked directly.  $\square$ \bigskip

Note that $\tau_5=1$ on $K=\mathbb{Q}(i)$ and $T(i)=T_2(i)=-i$, so the solution $(\xi, \eta)=(-i,i)$ of $\mathcal{C}_5$ is covered  by Theorem 5.1. \bigskip

If we substitute the expression in Theorem 5.1 for $\xi$ into the equation for $\mathcal{C}_5$ and simplify, we obtain:
$$(\eta^{4\tau_5}+2\eta^{3\tau_5}+4\eta^{2\tau_5}+3\eta^{\tau_5}+1)\eta^5=\eta^{\tau_5}(\eta^{4\tau_5}-3\eta^{3\tau_5}+4\eta^{2\tau_5}-2\eta^{\tau_5}+1).\eqno{(5.4)}$$
Thus, we have: \bigskip

\noindent {\bf Theorem 5.2.} {\it If 
$$g(X,Y)=(Y^4+2Y^3+4Y^2+3Y+1)X^5-Y(Y^4-3Y^3+4Y^2-2Y+1),$$
then $(X,Y)=(\eta, \eta^{\tau_5})$ is a point on the curve $g(X,Y)=0$.} \bigskip

From this we deduce the following.  \bigskip

\noindent {\bf Theorem 5.3.} {\it The roots of $p_d(x)$ are periodic points of the multi-valued algebraic function $\mathfrak{g}(z)$ defined by $g(z,\mathfrak{g}(z))=0$.  The period of $\eta$ with respect to the action of $\mathfrak{g}$ is the order of $\displaystyle \tau_5=\left(\frac{\mathbb{Q}(\eta)/K}{\wp_5} \right)$ in $\textrm{Gal}(\mathbb{Q}(\eta)/K)$.} \medskip

\noindent {\bf Remark.} See the Introduction of Part I for the definition of a periodic point of an algebraic function. \smallskip

\noindent {\it Proof.} Since $g(X,Y)$ has rational coefficients, applying $\tau_5^i$ ($1\le i \le n-1$) to the equation $g(\eta, \eta^{\tau_5}) = 0$ gives that
$$g(\eta, \eta^{\tau_5}) =g(\eta^{\tau_5}, \eta^{\tau_5^2}) = \dots = g(\eta^{\tau_5^{n-1}}, \eta)=0,$$
where $n = \textrm{ord}(\tau_5)$.  Thus, $\eta$ is one of the values of the iterate $\mathfrak{g}^{(n)}(\eta)$, i.e., is periodic with period $n$.  Any conjugate over $\mathbb{Q}$ of a periodic point of $\mathfrak{g}(z)$ is also a periodic point, and this proves the theorem.  $\square$ \bigskip

By Artin Reciprocity, the order of $\tau_5$  is equal to the order of $\wp_5$ in the quotient group $\textsf{A}/(\textsf{S}_{\wp_5'} \cap \textsf{P}_{f})$ (when $d \neq 4f^2$), where $\textsf{A}$ is the group of fractional ideals in $K$ which are relatively prime to $\wp_5'(f)$.  If this order is $n$, then there is an equation $\wp_5^n = (\frac{x+y\sqrt{-d}}{2})$, and since $y\sqrt{-d} \equiv x$ (mod $\wp_5'$), it follows that $\alpha = \frac{x+y\sqrt{-d}}{2} \equiv 2x/2 = x \equiv \pm 1$ (mod $\wp_5'$).  Therefore, when $d \neq 4f^2$, the period $n$ of the roots of $p_d(x)$ is the smallest positive integer $n$ for which there is an equation $4 \cdot 5^n =x^2+dy^2$ with $x \equiv \pm 1$ (mod $5$) and $(x,y) \mid 2$.  \medskip

The substitution $(X,Y) \rightarrow \left(\frac{-1}{X},\frac{-1}{Y}\right)$ represents an automorphism of the curve $g(X,Y)=0$, since
$$X^5 Y^5 g\left(\frac{-1}{X},\frac{-1}{Y}\right)=g(X,Y).$$
The equation connecting $t=X-\frac{1}{X}$ and $u=Y-\frac{1}{Y}$ in the function field of this curve is
$$h(t,u)=u^5-(6+5t+5t^3+t^5)u^4+(21+5t+5t^3+t^5)u^3-(56+30t+30t^3+6t^5)u^2$$
$$+(71+30t+30t^3+6t^5)u-120-55t-55t^3-11t^5;\eqno{(5.5)}$$
this follows from the calculation
$$-h(t,u)^2=\textrm{Res}_y(\textrm{Res}_x(g(x,y),x^2-tx-1),y^2-uy-1).$$
From $g(\eta,\eta^{\tau_5})=0$ and $\upsilon^{\tau_5}=\eta^{\tau_5}-\frac{1}{\eta^{\tau_5}}$ we obtain
$$h(\upsilon, \upsilon^{\tilde \tau_5}) = 0, \ \ \ \tilde \tau_5 = \tau_5 |_{\Omega_f} = \left(\frac{\Omega_f/\mathbb{Q}(\sqrt{-d})}{\wp_5}\right).$$
This yields the following result. \bigskip

\noindent {\bf Theorem 5.4.}  {\it If $\Omega_f$ is the ring class field of conductor $f$ (relatively prime to $5$) over the field $K=\mathbb{Q}(\sqrt{-d})$, where $-d=d_Kf^2$ and $\left(\frac{-d}{5}\right)=+1$, then $\Omega_f=K(\upsilon)$, where $\upsilon=\eta-\frac{1}{\eta}$ is a periodic point of the algebraic function $\mathfrak{f}(z)$ defined by $h(z,\mathfrak{f}(z))=0$, and $h(t,u)$ is given by equation (5.5).  The period of $\upsilon$ is the order of $\tilde \tau_5=\tau_5|_{\Omega_f}$ in $\textrm{Gal}(\Omega_f/K)$.}  \bigskip

Now we compare (5.4) with Ramanujan's modular equation
$$r^5(\tau)=r(5\tau)\frac{r^4(5\tau)-3r^3(5\tau)+4r^2(5\tau)-2r(5\tau)+1}{r^4(5\tau)+2r^3(5\tau)+4r^2(5\tau)+3r(5\tau)+1}$$
for $r(\tau)$.  Setting
$$\mathfrak{r}(z) = \frac{z(z^4-3z^3+4z^2-2z+1)}{z^4+2z^3+4z^2+3z+1},$$
we conclude from (5.4) and Theorem 4.5 that
$$\mathfrak{r}(\eta^{\tau_5})=\eta^5 = r^5(w/5) = \mathfrak{r}(r(w)). \eqno{(5.6)}$$
It is easily checked on Maple that the quintic extension of function fields $\mathbb{Q}(\zeta_5, z)/\mathbb{Q}(\zeta_5, \mathfrak{r}(z))$ is normal and cyclic, with generating automorphism
$$z \rightarrow \mathfrak{s}(z) = \frac{(\zeta+\zeta^2)z+1}{z+1+\zeta+\zeta^2},$$
where $\mathfrak{s}(z)=S^{-2}AS(z)=S^{-1}TS^{-1}(z)$ is an element of $G_{60}$.  It follows from (5.6) that
$$\eta^{\tau_5} = \mathfrak{s}^i(r(w)), \ \textrm{for some} \ \ i, \ 0 \le i \le 4.$$
From Corollary 4.7 and Theorem 4.8 we know that $i \neq 0$, since $\eta^{\tau_5} \in F_1$, but $r(w)$ generates $F$.  More specifically, we have the following.  \bigskip

\noindent {\bf Theorem 5.5.} {\it With notation as above, if $\xi=\zeta^j r(-1/w)$, $1 \le j \le 4$, we have the formula
$$r(w/5)^{\tau_5} = \mathfrak{s}^j(r(w))=T(\xi),$$
and $j$ is the unique integer (mod $5$) for which $\mathfrak{s}^j(r(w))$ is a root of $p_d(x)$.}  \medskip

\noindent {\it Proof.} We have that $\xi=\zeta^j r(-1/w)= S^j T(r(w))$, by the transformation formula for $r(-1/w)$, so
$T(\xi) = T S^j T (r(w))$.  On the other hand, $\mathfrak{s}(z) = S^{-1}TS^{-1}(z) = TST(z)$, since $(ST)^3 = 1$.  Therefore, $\mathfrak{s}^j(r(w)) = (TST)^j(r(w)) = TS^j T(r(w)) = T(\xi)$ since $T$ is its own inverse.  The above formula now follows from (5.2).  This proves that $\mathfrak{s}^j(r(w))$ is a root of $p_d(x)$, since $p_d(x)$ is stabilized by $T$.  There is only one value of $i$ for which $\mathfrak{s}^i(r(w))$ is a root of $p_d(x)$, since $T(\mathfrak{s}^i(r(w)))=S^i T(r(w)) = \zeta^i r(-1/w)$ must also be a root of $p_d(x)$.  $\square$  \bigskip
 
\noindent {\bf Remark.}  Since $\mathfrak{s}(z)=TST(z)$, $\mathfrak{s}(r(w))=TST(r(w))=TS(r(-1/w))=T(r(1-1/w))=r(-w/(w-1))$. Thus, $\mathfrak{s}^j(r(w))=r(w/(1-jw))$.  \bigskip

\noindent {\bf Example 1.} Consider Ramanujan's remarkable value
$$r(3i) = \sqrt{c^2+1}-c, \ \ 2c=\frac{60^{1/4}+2-\sqrt{3}+\sqrt{5}}{60^{1/4}-2+\sqrt{3}-\sqrt{5}} \sqrt{5}+1$$
established in  \cite{bch} and  \cite{bcz}, p.142.  A calculation on Maple shows that the minimal polynomial of $r(3i)=\zeta_5 r(4+3i)=\zeta r(w)$ is
\begin{align*}
m(x)&=x^{16}+38x^{15}-240x^{14}-300x^{13}-235x^{12}-726x^{11}+92x^{10}-1840x^9\\
&-675x^8+1840x^7+92x^6+726x^5-235x^4+300x^3-240x^2-38x+1,
\end{align*}
which is a factor of $G_{36}(x^5)$ in (3.11).  (Use the polynomial $H_{-36}(x)$ given in the proof of Proposition 3.2.)  Thus, $r(3i)$ is a linear fractional expression in some conjugate of $\eta=r\left(\frac{4+3i}{5}\right)$ with coefficients in $L=\mathbb{Q}(\zeta_5)$, and the minimal polynomial of the latter value is
$$p_{36}(x) = x^8 +x^6-6x^5+9x^4+6x^3+x^2+1,$$
from Table 1.  Using Maple to compare approximations of $r\left(\frac{4+3i}{5}\right)$ and the roots of $p_{36}(x)$, we find 
$$r\left(\frac{4+3i}{5}\right)=\frac{-i \omega^2}{2}+\frac{i\sqrt{3}}{2}-\frac{\omega}{4} \sqrt[4]{3}\left(\sqrt{4+2\sqrt{5}}+i\sqrt{-4+2\sqrt{5}}\right),\eqno{(5.7)}$$
with $\omega=\frac{-1+i\sqrt{3}}{2}$.  \smallskip

We determine the linear fractional expression in a root of $p_{36}(x)$ which will equal $r(3i)$.  Since
$$p_{36}(x) \equiv (x+3)^4(x^4+3x^3+x^2+2x+1) \ \ (\textrm{mod} \ 5),$$
the Frobenius automorphism $\tau_5$ has order $4$.  A calculation on Maple shows that
$$\mathfrak{s}^2(r(w))=\frac{(\zeta+\zeta^3)r(w)+1}{r(w)+1+\zeta+\zeta^3}=1.375418808...-(.899074105...)i$$
is the unique value $\mathfrak{s}^j(r(w))$ which is a root of $p_{36}(x)=0$. By Theorem 5.5 we have 
$$\eta^{\tau_5} = \mathfrak{s}^2(r(w))=\frac{(\zeta+\zeta^3)r(w)+1}{r(w)+1+\zeta+\zeta^3}=\frac{(1+\zeta^2)r(3i)+1}{\zeta^4 r(3i)+1+\zeta+\zeta^3}.$$
Inverting the linear fractional map in the last equality gives
$$r(3i)=\frac{(1+\zeta^3)\eta^{\tau_5}+\zeta}{\eta^{\tau_5}-\zeta-\zeta^3};$$
this is the desired expression for $r(3i)$.  Another calculation on Maple using (5.7) shows that
$$\eta^{\tau_5}=r\left(\frac{4+3i}{5}\right)^{\tau_5}=\frac{-i \omega}{2}-\frac{i\sqrt{3}}{2}+i\frac{\omega^2}{4} \sqrt[4]{3}\left(\sqrt{4+2\sqrt{5}}+i\sqrt{-4+2\sqrt{5}}\right).$$
This expresses $r(3i)$ in terms of $3$rd, $4$th, and $5$th roots of unity and shows that $\tau_5$ can be given by
$$\tau_5= \left(\sqrt[4]{3} \rightarrow -i\sqrt[4]{3}, i \rightarrow i, \sqrt{4+2\sqrt{5}} \rightarrow \sqrt{4+2\sqrt{5}}\right)|_{F_1}.$$ 
This proves formula (1.6) of the Introduction.  \bigskip

\noindent {\bf Remark.} In this example, $F=\Sigma_5 \Omega_{15}$ has degree $8h(-36)=16$ over $K = \mathbb{Q}(i)$, so its real subfield $F^+$ has degree $16$ over $\mathbb{Q}$ and the value $r(3i)$ generates $F^+$.  In particular, $K(r(3i)) = \Sigma_5 \Omega_{15}$.  Since $\sqrt{3} \in \Omega_3 \subset \Omega_{15}$ and $\sqrt{5} \in \Omega_5 \subset \Omega_{15}$, Ramanujan's formula shows that $60^{1/4} \in \Sigma_5 \Omega_{15}$.  On the other hand, $\Omega_3(60^{1/4})$ is a cyclic quartic extension of $\Omega_3$.  As in the proof of Theorem 4.6, there are only two cyclic quartic extensions of $\Omega_3$ contained in $\Sigma_5 \Omega_{15}$, namely, $\Sigma_5 \Omega_3 = \Omega_3(\zeta_5)$ and $\Omega_{15}$ (see Section 3); and the former is abelian over $\mathbb{Q}$.  Hence, we have $\Omega_{15}=K(\sqrt{3}, \sqrt[4]{60})$.  As a corollary, this shows that the rational primes which split completely in $\Omega_{15}$, which are the primes representable as $p=a^2+15^2 b^2$, are characterized by the two conditions $p \equiv 1$ (mod $12$) and $\left(\frac{60}{p}\right)_4=+1$.  $\square$  \medskip

Given that the period of $\eta$ in the above example is $n=4$, $p_{36}(x)$ can be calculated by a threefold iterated resultant, as in Part I, Section 3, pp. 727-730.  Namely, $p_{36}(x)$ is a factor of
$$R_4(x) = Res_{x_3}(Res_{x_2}(Res_{x_1}(g(x,x_1),g(x_1,x_2)),g(x_2,x_3)),g(x_3,x)).$$
Unfortunately, this calculation takes an extremely long time to complete, since $\textrm{deg}(R_4(x))=2\cdot5^4-1=1249$.  \medskip

To get around this difficulty, we let $g_1$ be the polynomial $g_1(X,Y) = Y^5 g(X,\frac{-1}{Y})$, i.e.,
$$g_1(X,Y)=  Y(Y^4-3Y^3+4Y^2-2Y+1)X^5+(Y^4+2Y^3+4Y^2+3Y+1).$$
The class number $h(-36)=2$, so $[F_1: K]=4$, hence $\textrm{Gal}(F_1/K) = \langle \tau_5 \rangle$, implying that $\tau_5^2=\rho$ on $F_1$.  Putting $\tau=\tau_5$, we have
$$g(\eta,\eta^\tau)=g(\eta^\tau,\eta^{\tau^2})=0.$$
However, $g(\eta^\tau,\eta^{\tau^2})=g(\eta^\tau,\eta^\rho)=g(\eta^\tau,-1/\eta)$, so that
$$g(\eta,\eta^\tau)=g_1(\eta^\tau,\eta)=0.$$
Therefore, $p_{36}(x)$ should be a factor of the resultant
\begin{align*}
\tilde R_2(x) &= Res_{x_1}(g(x,x_1),g_1(x_1,x))\\
& = -(x^2+1)(x^8+x^7+x^6-7x^5+12x^4+7x^3+x^2-x+1)\\
& \times (x^8+4x^7-x^6-14x^5+23x^4+14x^3-x^2-4x+1)\\
& \times (x^8-2x^7+x^6-4x^5+3x^4+4x^3+x^2+2x+1)\\
& \times (x^8+x^6-6x^5+9x^4+6x^3+x^2+1)\\
& \times (x^{16}+4x^{15}+29x^{12}-24x^{11}+86x^{10}-32x^9+105x^8\\
& \ \ +32x^7+86x^6+24x^5+29x^4-4x+1)\\
& = -(x^2+1)p_{51}(x) p_{91}(x) p_{24}(x) p_{36}(x) p_{96}(x).
\end{align*}
Hence, the discriminants with $d \in \{24, 36, 51, 91, 96\}$ are {\it all} the discriminants for which $\tau_5^2 = \rho$.  An analysis similar to the above for $d=36$ can be applied for these integers $d$ to yield formulas for the corresponding values of the Rogers-Ramanujan continued fraction $r(w)$, namely,
$$r(12+\sqrt{-6}), \ r\left(\frac{7+\sqrt{-51}}{2}\right), \ r\left(\frac{3+\sqrt{-91}}{2}\right), \ r(1+2\sqrt{-6}).$$
In addition, for small values of $n$, the $(n-1)$-fold iterated resultant
$$\tilde R_n(x) = Res_{x_{n-1}}(...(Res_{x_2}(Res_{x_1}(g(x,x_1),g(x_1,x_2)),g(x_2,x_3)),...,g_1(x_{n-1},x))$$
can be used to determine minimal polynomials of $r(w/5)$ for the values of $d \equiv \pm 1$ (mod $5$) for which $\rho \in \langle \tau_5 \rangle$ and $\tau_5^n =\rho$. \bigskip

\noindent {\bf Example 2.}  For example, $\tilde R_3(x)$ has degree $226$ and is the product of $(x^2+1)$ and $2$ factors of degree $4$, $3$ factors of degree $12$, $4$ factors of degree $24$, and one factor each of degree $36$ and $48$.  The degree $36$ factor is
\begin{align*}
p_{491}(x)&=x^{36}+28x^{35}+206x^{34}-324x^{33}+2163x^{32}+2080x^{31}+1600x^{30}\\
& +19440x^{29}+9145x^{28}+60876x^{27}+21486x^{26}-5532x^{25}+220279x^{24}\\
& +208904x^{23}+453304x^{22}-117152x^{21}-62271x^{20}+142940x^{19}\\
& +1116798x^{18}-142940x^{17}-62271x^{16}+117152x^{15}+453304x^{14}\\
& -208904x^{13}+220279x^{12}+5532x^{11}+21486x^{10}-60876x^9+9145x^8\\
& -19440x^7+1600x^6-2080x^5+2163x^4+324x^3+206x^2-28x+1,
\end{align*}
with discriminant $D=2^{316} 5^{153} 7^{16} 19^4 23^8 29^{16} 191^8 491^{18}$.  The value $d=491$ is a guess based on the conjecture at the end of Section 4.  This can be verified by factoring $p_{491}(x)$ modulo primes of the form $p=(x^2+491y^2)/4$, with $x+3y \equiv \pm 2$ (mod $5$) (assuming that $w=\frac{3+\sqrt{-491}}{2}$), to check that it splits into linear and quadratic factors.  For example, $p_{491}(x)$ factors into a product of linear polynomials modulo the primes $179=\frac{15^2+491}{4}, 3251=\frac{27^2+5^2 \cdot 491}{4}$, and $3989=45^2+2^2 \cdot 491$; while it splits into a product of $18$ linear factors and $9$ quadratics modulo $1237=\frac{23^2+3^2\cdot 491}{4}$, corresponding to the fact that $(\alpha)=\left(\frac{23+3\sqrt{-491}}{2}\right)$ satisfies $\alpha \equiv 1$, but $\alpha' \equiv 2$ (mod $\wp_5'$).  As an additional check, $\eta=r\left(\frac{3+\sqrt{-491}}{10}\right)$ is a root of $p_{491}(x)$ (to an accuracy of at least $60$ decimal places).  Note that $\textrm{ord}(\tau_5)=6$, since $\tau_5^3=\rho$ has order $2$, so the roots of $p_{491}(x)$ have period $6$ with respect to the action of $\mathfrak{g}(z)$.  This aligns with the fact that $4 \cdot 5^3 = 3^2 + 491$ and $4 \cdot 5^6=241^2+3^2 \cdot 491$ and that
$$\alpha_1=\frac{3+\sqrt{-491}}{2} \not \in \textsf{S}_{\wp_5'} \ \ \textrm{but} \ \ \alpha_2 = \frac{241+3\sqrt{-491}}{2} \in \textsf{S}_{\wp_5'}.$$

In general, it is more convenient to work with a lower degree polynomial derived from $p_d(x)$ using the fact that it is stabilized by the subgroup $H$.  First write $p_d(x)=x^{2h(-d)}t_d\left(x-1/x\right)$, which is possible since $p_d(x)$ is stabilized by $U(z)=-1/z$ (or $\eta^\rho=-1/\eta$ is an automorphism fixing $\Omega_f$).  Then $t_d(x)$ is a normal polynomial with root $\upsilon=\eta-1/\eta$ generating $\Omega_f$.  By (5.1), we can write $t_d(x-1) = x^{h(-d)} u_d\left(x+\frac{5}{x}\right)$.  This yields the polynomial $u_d(x)$ having degree $h(-d)$ and smaller discriminant.  In the above example we find 
\begin{align*}
u_{491}(x)= & x^9+10x^8-144x^7-840x^6+18354x^5-110972x^4+345800x^3\\
& -601496x^2+550293x-205102,
\end{align*}
whose discriminant is $D_1 = 2^{76} 7^2 29^4 191^2 491^4$.  It is straightforward to check that $7, 29, 191$ divide the index and $491$ does not (using Dedekind's method in \cite{ded}, pp. 214-218, for example), so we only have to exclude $q = 2$ and $q=29$ as divisors of $d$.  However, $h(-4 \cdot 29) = 6$ and $h(-491)=9$ yield that $d= 491 f^2$, where $f=2^a$.  If $a \ge 2$, then $h(-d)$ is even, while $h(-4 \cdot 491)= 27$, so the only possibility is $d=491$.  \medskip

A similar analysis was applied to check the polynomials in Tables 1 and 2.  \medskip

We will continue this discussion in Part III, by showing that the only irreducible factors of iterated resultants of the form $R_n(x)$ or $\tilde R_n(x)$ are the polynomials $x, x^2+1$, and $p_d(x)$, for $d \equiv \pm 1$ (mod $5$).  This will prove that the polynomial $p_{491}(x)$ given above actually is the minimal polynomial of $r(w/5)$ for $w=\frac{3+\sqrt{-491}}{2}$.

\medskip

\noindent Dept. of Mathematical Sciences, LD 270

\noindent Indiana University - Purdue University at Indianapolis (IUPUI)

\noindent Indianapolis, IN 46202

\noindent {\it e-mail: pmorton@iupui.edu}

\end{document}